%% file: article.tex
\documentclass[11pt,twoside]{article}

\input{shared}

\usepackage[font=footnotesize]{caption} 

\begin{document}

\date{}
\maketitle
\thispagestyle{empty}
\setlength{\headheight}{13.6pt}
\pagestyle{fancy}
\fancyhead{}

\fancyhead[OC]{Fast and Scalable FFT-Based GPU-Accelerated Algorithms for Hessian Actions}
\fancyhead[EC]{Sreeram Venkat, Milinda Fernando, Stefan Henneking, and Omar Ghattas}

\begin{abstract}
We present an efficient and scalable algorithm for performing
matrix-vector multiplications ("matvecs") for block Toeplitz
matrices. Such matrices, which are shift-invariant with respect to
their blocks, arise in the context of solving inverse problems
governed by autonomous systems, and time-invariant systems in
particular. In this article, we consider inverse problems that infer
unknown parameters from observational data of a linear time-invariant
dynamical system given in the form of partial differential equations
(PDEs). Matrix-free Newton-conjugate-gradient methods are often the
gold standard for solving these inverse problems, but they require
numerous actions of the Hessian on a vector. Matrix-free adjoint-based
Hessian matvecs require solution of a pair of linearized
forward/adjoint PDE solves per Hessian action, which may be
prohibitive for large-scale inverse problems. Time invariance of the
forward PDE problem leads to a block Toeplitz structure of the
discretized parameter-to-observable (p2o) map defining the mapping
from inputs (parameters) to outputs (observables) of the PDEs. This
block Toeplitz structure enables us to exploit two key properties: (1)
compact storage of the p2o map and its adjoint; and (2) efficient fast
Fourier transform (FFT)-based Hessian matvecs. The proposed algorithm
is mapped onto large multi-GPU clusters and achieves more than 80
percent of peak bandwidth on NVIDIA A100 GPUs. Excellent weak scaling
is shown for up to 48 A100 GPUs. For the targeted problems, the
implementation executes Hessian matvecs within fractions of a second,
which is orders of magnitude faster than can be achieved by 
conventional matrix-free Hessian matvecs via forward/adjoint PDE
solves. 

\end{abstract}

\newpage

\section{Introduction}\label{sec:intro}

\input{introduction}

\section{Motivation and Background}\label{sec:background}

\input{background}

\section{Methods}\label{sec:methods}

\input{methods}

\section{Numerical Results}\label{sec:results}

\input{numerical_results}

\section{Conclusions}\label{sec:conclusions}

\input{conclusions}

\section*{Acknowledgments}
This research was supported in part by DOE ASCR grants DE-FOA-0002704
and DE-SC0021239 and DOD MURI grant FA9550-21-1-0084. Supercomputing
resources were provided by the Texas Advanced Computing Center (TACC)
at UT Austin on its \emph{Frontera} and \emph{Lonestar6} systems. This material is
based upon work supported by the National Science Foundation Graduate
Research Fellowship under Grant No.~DGE 2137420. We would also like to
thank Professor Mike Giles of the University of Oxford for the
insightful discussion that led to the local matvec algorithm presented
in~\cref{sec:LocalMatvec}. 

\appendix
\crefalias{section}{appendix}
\section{Appendix: Alternate Algorithm for Local Matvecs}\label{sec:appendix-alt-matvec}
\input{appendix}

\bibliographystyle{siamplain}

\end{document}

%% file: shared.tex
\input{preamble}

\title{
Fast And Scalable FFT-Based GPU-Accelerated Algorithms for Block-Triangular Toeplitz Matrices With Application to Linear Inverse Problems Governed by Autonomous Dynamical Systems
}

\author[1]{Sreeram Venkat\thanks{Corresponding author: srvenkat@utexas.edu}}
\author[1]{Milinda Fernando}
\author[1]{Stefan Henneking}
\author[1,2]{Omar Ghattas}
\affil[1]{Oden Institute, The University of Texas at Austin}
\affil[2]{Walker Department of Mechanical Engineering, The University of Texas at Austin}

%% file: preamble.tex
\usepackage[T1]{fontenc}
\usepackage[utf8]{inputenc}
\usepackage[english]{babel}
\usepackage[margin=1in]{geometry}

\usepackage{amsfonts}
\usepackage{amsmath}
\usepackage{fancyhdr}
\usepackage{graphicx}
\usepackage{xcolor}
\usepackage{authblk}
\usepackage{adjustbox}
\usepackage{setspace}
\usepackage{multirow}
\usepackage{multicol}
\usepackage{epstopdf}
\usepackage{url}
\usepackage{algorithm}
\usepackage[noend]{algpseudocode}
\usepackage{caption}
\usepackage{subcaption}
\usepackage{float}
\usepackage{placeins}
\usepackage{bm}
\usepackage{booktabs}
\usepackage{physics}
\usepackage{cancel}
\usepackage{amsopn}

\usepackage[pdftex,hidelinks]{hyperref}
\hypersetup{colorlinks=true, citecolor=blue, linkcolor=blue, urlcolor=blue}

\usepackage[capitalize]{cleveref}
\usepackage{capt-of}
\usepackage[normalem]{ulem}

\usepackage{pgfplots}
\DeclareUnicodeCharacter{2212}{−}
\usepgfplotslibrary{groupplots,dateplot}
\usetikzlibrary{patterns,shapes.arrows}
\pgfplotsset{compat=newest}
\pgfplotsset{compat=newest,
    /pgfplots/ybar legend/.style={
    /pgfplots/legend image code/.code={%
    \draw[##1,/tikz/.cd,bar width=3pt,yshift=-0.2em,bar shift=0pt]
    plot coordinates {(0cm,0.8em)};},
   },
}
\pgfplotsset{hide scale/.style={
/pgfplots/xtick scale label code/.code={},
/pgfplots/ytick scale label code/.code={}}}

\ifpdf{}
  \DeclareGraphicsExtensions{.eps,.pdf,.png,.jpg}
\else
  \DeclareGraphicsExtensions{.eps}
\fi
\usetikzlibrary{patterns}

\DeclareMathOperator*{\argmin}{arg\,min}

\usepackage{amsthm}
\theoremstyle{definition}
\newtheorem{remark}{Remark}

\newcommand{\R}{\mathbb{R}}

\newcommand{\N}{\mathbb{N}}

\newcommand{\blocktoeptxt}{F}
\newcommand{\blocktoep}{\vb{\blocktoeptxt}}
\newcommand{\paramtxt}{m}
\newcommand{\datatxt}{d}
\newcommand{\paramvec}{\vb{\paramtxt}}
\newcommand{\datavec}{\vb{\datatxt}}
\newcommand{\obsvec}{\vb{\datatxt}^{\text{obs}}}

\newcommand{\toepblocktxt}{\widetilde{\blocktoeptxt}}
\newcommand{\toepblock}{\widetilde{\vb{\blocktoeptxt}}}
\newcommand{\tparamvec}{\widetilde{\vb{\paramtxt}}}
\newcommand{\tdatavec}{\widetilde{\vb{\datatxt}}}
\newcommand{\tparamtxt}{\widetilde{\paramtxt}}
\newcommand{\tdatatxt}{\widetilde{\datatxt}}

\newcommand{\numgrid}{N_g}
\newcommand{\numstate}{N_u}
\newcommand{\numtime}{N_t}
\newcommand{\numparam}{N_m}
\newcommand{\numdata}{N_d}

\newcommand{\numprocrows}{r}
\newcommand{\numproccols}{c}
\newcommand{\procrowblock}{n_d}
\newcommand{\proccolblock}{n_m}
\newcommand{\proc}{p}

\definecolor{RED}{RGB}{190, 30, 49}
\definecolor{GREEN}{RGB}{0, 171, 86}
\definecolor{BLUE}{RGB}{11, 78, 179}

\newcommand{\revstart}{\color{black}}
\newcommand{\revend}{\color{black}}
\newcommand{\rev}[1]{\textcolor{black}{#1}}

\newcommand{\be}{\begin{equation}}
\newcommand{\ee}{\end{equation}}

\newcommand{\eq}[1]{(\ref{eq:#1})}

\newcommand{\bb}[1]{\mathbb{#1}}
\newcommand{\mc}[1]{\mathcal{#1}}

\newcommand{\p}{\partial}

\crefname{subsection}{Section}{Sections}

%% file: introduction.tex
\revstart
Inverse problems are solved widely throughout computational science to determine unknown parameters or inputs to models of physical phenomena from observations~\cite{engl1996regularization,isakov2006inverse,vogel2002computational,kirsch2011introduction,banks2012estimation}. In many inverse problems of scientific interest, the observations are often sparse and noisy while the parameter or input fields to be determined are high-dimensional (usually discretizations of continuous fields). As a result, these inverse problems are frequently ill-posed, and  specialized mathematical and computational techniques are required to treat them properly~\cite{stuart2010inverse, kaipio2006statistical, tarantola2005inverse}. We refer readers to~\cite{ghattas2021learning} (Part 2) for an overview of state-of-the-art computational methods for large-scale inverse problems. 
\revend
Adjoint-based matrix-free Newton-conjugate-gradient methods are often the gold standard for solution of \rev{such} inverse problems~\cite{ghattas2021learning, villa2021hippylib}.
However, they typically require numerous actions of the Hessian matrix on a vector, each of which amounts to solution of a pair of forward and adjoint problems.
For inverse problems governed by partial differential equations (PDEs) with high-rank Hessians, the resulting number of
forward/adjoint PDE solves may be computationally prohibitive. 
\revstart
Constructing reduced-order models for these problems can also be difficult due to the high dimensionality of the parameter fields and the Kolmogorov $N$-width problem~\cite{greif2019kolmogorov}, which asserts the nonexistence of a low-dimensional linear subspace embedding.
\revend

In recent years, several methods that address high-rank Hessians have been developed, including those that exploit the pseudo-differential~\cite{nammour2011phd, demanet2012probing}, augmented Lagrangian~\cite{alger2017augmented}, product-convolution~\cite{alger2019product}, H-matrix~\cite{ambartsumyan2020hierarchical}, and point spread function~\cite{alger2024point} structure of the Hessians of particular classes of inverse problems.

In this paper, we show that the Hessian structure can be particularly well-exploited for Hessians governed by \emph{autonomous systems}. {The evolution of such systems with respect to any given input may depend on the system's current state but does not explicitly depend on the independent variable.} Autonomous systems can arise, for example, in the context of inverse problems for time-invariant dynamical systems, which are a subclass of autonomous systems where the independent variable is time~\cite{kloeden2011autonomous, gavin2018notes, henneking2025goal}. The autonomous system structure then translates into a \emph{shift invariance} of the corresponding discrete system. In particular, with both parameter field and observables defined in space-time, the discrete parameter-to-observable (p2o) map and its adjoint exhibit shift invariance with respect to the time-stepping. The corresponding matrices are \emph{block Toeplitz}. For time-invariant systems, causality additionally implies the p2o map and its adjoint are lower- and upper-triangular block Toeplitz, respectively. Recognizing this structure enables two {properties}: 1) compact representation of the p2o map and its adjoint, and 2) fast application of the Hessian {via} scalable multi-GPU fast Fourier transform (FFT){-accelerated} matvecs. Compact representation follows directly from the definition of block Toeplitz matrices. Fast {matvecs are} achieved by embedding the block Toeplitz matrix within a block circulant matrix, which is diagonalized by the discrete Fourier transform (DFT). The {matvec} then becomes an elementwise vector operation in Fourier space.

{We show that these FFT-based Hessian matvecs can be implemented
efficiently} {on multi-GPU clusters}. Moreover, because the FFT is a unitary operator, the action of the adjoint p2o map corresponds to simply applying the complex conjugate in Fourier space, eliminating the need to separately store the Fourier-transformed forward and adjoint maps. Exploiting the block-triangular Toeplitz structure in this way yields memory savings proportional to the number of time steps $\numtime$ and a computational speedup of $\mathcal{O}(\numtime / \log \numtime)$. In the context of explicit methods {for time-dependent differential equations}, the number of time steps is often very large due to the Courant--Friedrichs--Lewy (CFL) condition, making the savings of the algorithm substantial.

While the classical FFT algorithm for matvecs involving Toeplitz matrices is well known~\cite{gray2006toeplitz}, similar algorithms for {general block Toeplitz matrices---i.e.~block Toeplitz matrices where the blocks themselves do not have any special structure---have yet to be established}. Many authors have extended the FFT-based Toeplitz matvec algorithm to block Toeplitz matrices where each block is itself a Toeplitz matrix~\cite{lee1986fast, barrowes2001fast, kazeev2013multilevel, hasan2012multiway, yagle2001fast, wax1983efficient}. These ``multilevel-'' or ``recursive-'' Toeplitz matrices arise in scattering problems and optimal surface interpolation~\cite{lee1986fast, barrowes2001fast}. In~\cite{gallivan1996high}, the authors discuss extensions of factoring algorithms to symmetric positive definite block Toeplitz matrices; however, the block Toeplitz matrices found in inverse problems are usually not even square matrices. The use of GPUs to accelerate computations involving Toeplitz matrices (through FFTs) is also well documented~\cite{luo2021optimization}, though the extension to {general block Toeplitz matrices} and {the mapping onto multi-GPU clusters} to perform large-scale computations is yet to be {reported}.

{
The main contributions of this paper are: (1) novel algorithms for
efficient GPU-accelerated FFT-based matvecs of general block Toeplitz
matrices as arise in linear inverse problems governed by autonomous
dynamical systems; (2) extension of the algorithms to a multi-GPU
framework; (3) detailed roofline performance analysis for the
algorithms; (4) strong and weak scalability study for up to 48 GPUs;
and (5) complexity analysis of the algorithms for single-GPU and
multi-GPU execution.  The contributions made in this paper enable
efficient and scalable FFT-based Hessian matvecs for solving
large-scale inverse problems governed by autonomous dynamical systems.
To our knowledge, this is the first paper to recognize and
exploit the block Toeplitz structure of Hessian operators arising in
inverse problems governed by time-invariant dynamics via fast
algorithms and their extension to the multi-GPU setting.  }

%% file: background.tex
Our motivation comes from the need to perform Hessian matvecs to solve large-scale inverse problems governed by autonomous systems---time-invariant systems in particular. In this paper, we are concerned with systems governed by PDEs. For such problems, we briefly review the structure of the inverse problem and discuss the case of systems with time invariance that give rise to shift-invariant discrete operators~\cite{henneking2025goal}. Before presenting our algorithm in~\Cref{sec:methods}, we outline the classical FFT-based matvec algorithm for general Toeplitz matrices.

\revstart
\begin{remark}
    Throughout the paper, boldface text will be used to refer to discrete matrices and vectors (including block matrices and vectors). Individual elements of matrices or vectors will be denoted with regular text. 
\end{remark}
\revend

\subsection{Linear Time-Invariant Dynamical System}
\label{sec:LinearTimeInvariantDynamicalSystem}

While our algorithm applies to linear autonomous dynamical systems in general, we consider for illustrative purposes the case of a linear time-invariant (LTI) dynamical system of the form
\be
\begin{aligned}
	\frac{\partial u}{\partial t} &= \mc A u + \mc C m && \text{in } \Omega \times (0,T), \\
	u &= u_0 && \text{in } \Omega \times \{ 0 \}, \\
	d &= \mc B u && \text{in } \Omega \times (0,T),
\end{aligned}
\label{eq:LTI}
\ee
with appropriate boundary conditions on $\p \Omega \times (0,T)$, where $\Omega$ is the spatial domain, $(0,T)$ is the time domain, $u(x,t)$ is the system's state with initial state $u_0(x)$, parameter (input) $m(x,t)$ represents the source or forcing of the system and is independent of the state, and both $\mc A$ and $\mc C$ are time-invariant differential operators; $d(x,t)$ describes the observables (output) of the system, which are extracted from the state $u$ via a time-invariant observation operator $\mc B$.
\revstart
Examples of such LTI systems include heat transfer, diffusion, porous media flow, and wave propagation where $m$ represents a source term. In~\cite{henneking2025bell}, the algorithms developed in this paper are applied to an LTI system of coupled acoustic--gravity wave equations, where the parameter $m$ represents a boundary source for the velocity field.
\revend

Consider a discrete version of the LTI system obtained by discretizing \eq{LTI} in time with a single-step explicit method,
\be
	\vb{u}_{k+1} = \vb{A} \vb{u}_k + \vb{C} \vb{m}_k, \quad k=0,1,\cdots,\numtime-1,
\label{eq:LTI-time-stepping}
\ee
where $\vb{u}_k \in \bb{R}^{\numstate}$, $\vb{m}_k \in \bb{R}^{\numparam}$, and the discrete time-stepping operator $\vb{A} \in \bb{R}^{\numstate \times \numstate}$ and $\vb{C} \in \bb{R}^{\numstate \times \numparam}$ both depend on the particular time-stepping scheme.\footnote{For example, forward Euler time-stepping implies $\vb{A}$ and $\vb{C}$ respectively spatially discretize $(\mc I + \Delta t \mc A)$ and $\Delta t \mc C$, where $\mc I$ is identity and $\Delta t$ is the (uniform) time step size. \rev{For backward Euler time-stepping, $\vb{A}$ and $\vb{C}$ respectively spatially discretize $(\mc I - \Delta t \mc A)^{-1}$ and $\Delta t (\mc I - \Delta t \mc A)^{-1}\mc C$.}} Then, using \eq{LTI-time-stepping} the discretized LTI system can be written as follows:

\be
\label{eq:LTI-time-discrete}
\begin{split}
	\vb{u}_1 &= \vb{A} \vb{u}_0 + \vb{C} \vb{m}_0 , \\
	\vb{u}_2 &= \vb{A} \vb{u}_1 + \vb{C} \vb{m}_1 , \\
		&= \vb{A} (\vb{A} \vb{u}_0 + \vb{C} \vb{m}_0) + \vb{C} \vb{m}_1 = \vb{A}^2 \vb{u}_0 + \vb{A}^1 \vb{C} \vb{m}_0 + \vb{A}^0 \vb{C} \vb{m}_1 , \\[-2pt]
		\vdots \\[-10pt]
	\vb{u}_{k+1} &= \vb{A}^{k+1} \vb{u}_0 + \sum_{i=0}^k \vb{A}^i \vb{C} \vb{m}_{k-i} , \\[-3pt]
	\vb{d}_{k+1} &= \vb{B} \vb{u}_{k+1} ,
\end{split}
\ee
where $\vb{d}_k \in \bb{R}^{\numdata}$ and $\vb{B} \in \bb{R}^{\numdata \times \numstate}$ is the discrete observation operator.

Without loss of generality, assume homogeneous initial condition $\vb{u}_0 = 0$. We can then write the discretized LTI system in the following way:

\be
	\left[ 
	\begin{array}{@{\hskip 2pt}c@{\hskip 2pt}}
		\vb{d}_1 \\
		\vb{d}_2 \\
		\vdots \\
		\vb{d}_{k+1} \\
		\vdots \\
		\vb{d}_{\numtime}
	\end{array}
	\right]
	=
	\left[ 
	\begin{array}{@{\hskip 2pt}cc@{\hskip 4pt}c@{\hskip 4pt}c@{\hskip 4pt}c@{\hskip 4pt}c@{\hskip 2pt}}
		\vb{B} \vb{A}^0 \vb{C}  \\
		\vb{B} \vb{A}^1 \vb{C} & \vb{B} \vb{A}^0 \vb{C}  \\
		\vdots & \vdots & \ddots  \\
		\vb{B} \vb{A}^k \vb{C} & \vb{B} \vb{A}^{k-1} \vb{C} & \cdots & \vb{B} \vb{A}^0 \vb{C} \\
		\vdots & \vdots & & \vdots & \ddots  \\
		\vb{B} \vb{A}^{\numtime-1} \vb{C} & \vb{B} \vb{A}^{\numtime-2} \vb{C} & \cdots & \vb{B} \vb{A}^{\numtime-(k+1)} \vb{C} & \cdots & \vb{B} \vb{A}^0 \vb{C}
	\end{array}
	\right]
	\left[ 
	\begin{array}{@{\hskip 2pt}c@{\hskip 2pt}}
		\vb{m}_0 \\
		\vb{m}_1 \\
		\vdots \\
		\vb{m}_k \\
		\vdots \\
		\vb{m}_{\numtime-1}
	\end{array}
	\right] .
\label{eq:LTI-matrix}
\ee

We define $\vb{F}_{ij} := \vb{B} \vb{A}^{i-j} \vb{C} \in \bb{R}^{\numdata \times \numparam}$, $i,j=1,2,\ldots,\numtime, i \ge j$; the LTI system \Cref{eq:LTI-time-stepping} can then be written more compactly as follows:

\begin{equation}
	\left[ \begin{array}{c}
	\vb{d}_1 \\[2pt]
	\vb{d}_2 \\[6pt]
	\vb{d}_3 \\[1pt]
	\vdots \\[3pt]
	\vb{d}_{\numtime}
	\end{array} \right]
	=
	\left[ \begin{array}{ccccc}
	\vb{F}_{11} & 0 & 0 & \cdots & 0 \\[2pt]
	\vb{F}_{21} & \vb{F}_{11} & 0 & \cdots & 0 \\
	\vb{F}_{31} & \vb{F}_{21} & \vb{F}_{11} & \ddots & \vdots \\
	\vdots & \vdots & \ddots & \ddots & 0 \\[2pt]
	\vb{F}_{\numtime,1} & \vb{F}_{\numtime-1,1} & \cdots & \vb{F}_{21} & \vb{F}_{11}
	\end{array} \right]
	\hskip 5pt
	\left[ \begin{array}{c}
	\vb{m}_0 \\[2pt]
	\vb{m}_1 \\[6pt]
	\vb{m}_2 \\[1pt]
	\vdots \\[3pt]
	\vb{m}_{\numtime-1}
	\end{array} \right] ,
\label{eq:ShiftInvariance}
\end{equation}
or very concisely as
\be
	\datavec := \blocktoep \paramvec .
\label{eq:p2o}
\ee

We refer to $\blocktoep$ as the (discrete) \emph{parameter-to-observable} (p2o) map; $\paramvec$ is the parameter vector and $\datavec$ is the vector of observables or data vector.

Then,

\begin{itemize}
	\item $\paramvec \in \bb R^{\numparam \numtime}$ with blocks $\paramtxt_j \in \bb R^{\numparam}$, $j = 1,2, \ldots, \numtime$;
	\item $\datavec \in \bb R^{\numdata \numtime}$ with blocks $\datatxt_i \in \bb R^{\numdata}$, $i = 1,2, \ldots, \numtime$;
	\item $\blocktoep \in \bb R^{(\numdata \numtime) \times (\numparam \numtime)}$ with blocks $\blocktoeptxt_{ij} \in \bb{R}^{\numdata \times \numparam}$, $i,j = 1,2, \ldots, \numtime$.

\end{itemize}

It is clear from \Cref{eq:ShiftInvariance} that the p2o map $\blocktoep$ is shift-invariant with respect to its blocks $\blocktoeptxt_{ij}$. In particular, $\blocktoep$ is \emph{block Toeplitz}. Additionally, time causality implies that $\blocktoep$ is block lower-triangular. The next two sections describe how this special structure of $\blocktoep$ can be effectively exploited in the context of solving inverse problems.

\revstart
Note that the methodology is easily extended to multi-step explicit and implicit methods. In the case of multi-step methods, the operator $\vb{A}$ will represent the application of the spatial PDE operator over the intermediate time steps. In the case of implicit methods, the time stepping operator A requires solution of a linear system.  Moreover, $\vb{u}_k$ (and thus $\vb{d}_k$) may depend on $\vb{m}_k$; by taking the observation vector to be $[\vb{d}_0,\vb{d}_1,\dots,\vb{d}_{N_t}]^T$ and the parameter vector to be $[\vb{m}_0,\vb{m}_1,\dots,\vb{m}_{N_t}]^T$, the block-triangular Toeplitz structure of $\blocktoep$ is recovered.
\revend

\subsection{Hessian Matvec for Inverse Problems}
\label{sec:HessianMatvecforInverseProblems}

Given a system of the form \eq{p2o}, we consider the inverse problem of inferring the parameters $\paramvec$ from the observed data $\obsvec$.
This inverse problem can be solved by casting it as a quadratic optimization problem of minimizing the regularized data misfit:

\begin{equation}\label{eq:ObjectiveFunctional}
	\min_{\paramvec} J(\paramvec) := 
	\frac{1}{2} \| \blocktoep \paramvec - \obsvec  \|^2
	+ \frac{\alpha}{2} \| \paramvec \|_{\vb R}^2 ,
\end{equation}
where $\alpha>0$, the regularization operator $\vb R \in \bb{R}^{(\numparam \numtime) \times (\numparam \numtime)}$ is positive definite, and \rev{$\|\cdot\|_{\vb{R}}$ denotes the weighted $L^2$ norm.}

The goal is to minimize the objective $J(\paramvec)$ in \eq{ObjectiveFunctional} which, through p2o map $\blocktoep$, is constrained by the LTI system \eq{LTI-time-stepping}. Minimization of \eq{ObjectiveFunctional} defines a linear inverse problem:
\begin{equation}\label{eq:LinearInverseProblem}
	(\blocktoep^*  \blocktoep + \alpha \vb R) \paramvec = \blocktoep^*  \obsvec ,
\end{equation}
where $\blocktoep^*$ denotes the adjoint of the p2o map,\footnote{Note that analogous to the p2o map $\blocktoep$, the adjoint p2o map $\blocktoep^*$ is block Toeplitz; however, $\blocktoep^*$ is block upper-triangular. Moreover, applying $\blocktoep^*$ to a vector involves the adjoint operator $\vb{A}^*$ of the governing LTI system.} and $\vb H := \blocktoep^*  \blocktoep + \alpha \vb R \in \bb{R}^{(\numparam \numtime) \times (\numparam \numtime)}$ is the Hessian.

If the action of the Hessian on a vector is available, the inverse problem \eq{LinearInverseProblem} can be solved efficiently by iterative methods such as conjugate gradients (CG).
When preconditioned by the regularization operator, the number of iterations typically scales with the effective rank of the preconditioned Hessian of the data misfit, $\alpha^{-1} \vb R^{-1} \blocktoep^* \blocktoep$, since the resulting operator has the structure of a compact perturbation of the identity~\cite{ghattas2021learning}.

For the iterative solution to be practical, it is paramount that the Hessian matvec can be carried out efficiently at each iteration.

For large-scale inverse problems, the action of the Hessian on a vector is typically formed in a matrix-free way, i.e.~$\vb H$ is never explicitly constructed, since doing so requires either $\numparam \numtime$ forward solves or $\numdata \numtime$ adjoint solves, whichever is less expensive.
However, each Hessian matvec comes at the cost of a pair of forward and adjoint solves of the governing PDE system,
which can be very expensive when performed repeatedly.

For the LTI system \eq{LTI-time-stepping}, each forward solve requires $\numtime$ applications of the time-stepping forward operator $\vb{A}$.
Analogously, each adjoint solve requires $\numtime$ applications of a time-stepping adjoint operator $\vb{A}^*$.\footnote{For the LTI system, the cost of applying the adjoint operator is similar to the cost of applying the forward operator.}
This cost can make solving the inverse problem prohibitively expensive.

For example, consider inverse problem \eq{LinearInverseProblem} of size $\numdata = \mc O(10^2)$, $\numparam = \mc O(10^{6})$, $\numtime = \mc O(10^4)$.\footnote{Large-scale inverse problems for complex physical systems are often data-sparse---i.e.~the number of observers (sensors) is limited---while having high-dimensional parameter fields $(\numparam \gg \numdata)$.} Unless the Hessian has low rank, solving an inverse problem of this size with traditional methods is extremely challenging (indeed, it may require up to $\mc O(10^6)$ Hessian matvecs, i.e.~up to $\mc O(10^{10})$ applications of the time-stepping forward and adjoint operators $\vb{A}$ and $\vb{A}^*$).\footnote{\rev{The maximum effective rank of the Hessian is $N_dN_t \sim 10^6$. Using PDE-based Hessian actions, each action requires 2 PDE solves; time-stepping with $N_t \sim 10^4$ leads to $\mc O(10^{10})$ applications of $\vb{A}$.}}

To make this notion more concrete and be able to better compare the computational cost of the conventional approach---the matrix-free Hessian matvec via a pair of forward/adjoint PDE solves---to the proposed method, we estimate the number of floating-point operations (FLOPs) needed to apply the p2o map for a particular example.
For estimating FLOPs, we assume values for the number of parameters, time steps, grid points, and other variables, that represent practical large-scale problems including our own target applications.
Consider elastic wave propagation in displacement form
(e.g.~\cite{bedford2023elastic}) discretized with a 27-point stencil
in a uniformly refined three-dimensional cube with $\numgrid$ grid
points.\footnote{Discretization with a 27-point stencil is comparable
to a trilinear hexahedral finite element discretization. Since
higher-order discretization makes applying the discretized PDE
operator even more expensive, the 27-point stencil can be used as a lower bound for the speedup of our method.}
The spatial state dimension is $\numstate = 3 \numgrid$ (3 degrees of freedom (DOFs) per grid point).
We assume that the spatiotemporal parameter field is a (scalar-valued) field spatially discretized on the top surface of the spatial domain, $\numparam = \numgrid^{2/3}$, and observations are taken at a small number of sensors, $\numdata \ll \numparam$.
The cost for applying the spatially discretized PDE operator once is approximately $81 \numstate$~FLOPs.\footnote{Using a 27-point stencil with 3 state DOFs per grid point, each state DOF is connected to (at most) 81 DOFs.}
With a classic explicit RK4 time-stepping scheme, the spatially discretized operator is applied four times per time step, costing $324 \numstate$~FLOPs per time step.
Each forward (or adjoint) PDE solve involves $\numtime$ time steps, so the total cost per PDE solve is approximately $324 \numstate \numtime$~FLOPs.
Assuming the number of grid points is $\numgrid = 10^9$ and the number of time steps is $\numtime = 10^4$, each PDE solve costs $9.72 \cdot 10^{15}$~FLOPs or 9.72~petaFLOPs.
To estimate the number of PDE solves needed for solving the inverse
problem, we have to estimate the effective rank $r$ of the
(preconditioned) Hessian
\revstart 
which determines the approximate number of iterations for an iterative linear solver such as CG to converge on~\cref{eq:LinearInverseProblem} when regularization preconditioning is used~\cite{ghattas2021learning}.
\revend
Assuming $\numdata = 10^2$ sensors and that
10\% of the data are informative, $r = \numdata \numtime / 10 = 10^5$.
Then, solving the inverse problem in the conventional way costs
$2r \cdot 9.72 \cdot 10^{15} \approx 1.944 \cdot 10^{21}$~FLOPs or
1944~exaFLOPs.\footnote{Note that this analysis neglected the cost of
applying (or preconditioning with) the regularization operator $\vb
R$, which amounts to an elliptic solve and is typically significantly
cheaper than solving the forward or adjoint PDE system.}

In the next section, we propose a computationally much more efficient
approach for Hessian matvecs that makes solving inverse problems
governed by autonomous dynamical systems {\em orders of magnitude
cheaper} at this scale.

\subsection{Inverse Problems Involving Shift-Invariant Systems}\label{sec:InverseProblemsInvolvingShiftInvariantSystems}

Our approach exploits the shift invariance of \eq{ShiftInvariance} and trades computer storage for computational efficiency.
Assuming the above estimates for $\numdata$, $\numparam$, and $\numtime$, then formally the Hessian $\vb H$ and the p2o map $\blocktoep$ each have $\mc O(10^{10})$ columns.
Clearly, pre-computing and storing these dense matrices naively column-by-column by performing $\mc O(10^{10})$ matvecs is not a feasible option.
However, recognizing the block Toeplitz structure of the p2o map enables two key properties: (1) compact storage of the p2o map and its adjoint; and (2) efficient FFT-based Hessian matvecs. In particular:

\begin{itemize}
	\item $\blocktoep$ can be pre-computed by only $\numparam$ forward solves to obtain the first block column or $\numdata$ adjoint solves to obtain the last block row. \rev{In comparison, the naive approach would require performing $\numparam \numtime$ forward solves to compute $\blocktoep$ column-by-column or $\numdata \numtime$ adjoint solves to compute $\blocktoep^*$ column-by-column};
	\item $\blocktoep$ can be compactly stored in $\mc{O}(\numparam \numdata \numtime)$ memory, a savings of $\mc{O}(\numtime)$;
	\item $\blocktoep$ can be efficiently applied to a vector by a specialized FFT-based matvec algorithm at the cost of $\mc{O}(\numparam \numdata \numtime \log \numtime)$, a speedup of $\mc{O}(\numtime / \log \numtime)$ over \revstart a dense matvec\revend;
	\item $\blocktoep^*$ does not need to be computed or stored separately because the same FFT-based matvec algorithm with only minor modifications efficiently applies $\blocktoep^*$ to a vector, reducing the cost of the adjoint p2o matvec to $\mc{O}(\numparam \numdata \numtime \log \numtime)$.
\end{itemize}

Once the shift-invariant p2o map has been pre-computed and stored
compactly, the efficiency of the FFT-accelerated $\blocktoep$ and
$\blocktoep^*$ matvecs implies that the Hessian matvec becomes
relatively cheap, because it no longer involves applying the forward
or adjoint PDE time-stepping operators.
The computational cost of pre-computing $\blocktoep$ and $\blocktoep^*$ is small compared to the cost of solving the inverse problem using the conventional Hessian matvec method, as we show below.

Using the elastic wave propagation example and estimates
from \cref{sec:HessianMatvecforInverseProblems}, the one-time cost of
pre-computing the p2o operator is $\numdata$ adjoint solves or
$\numdata \cdot 9.72 \cdot 10^{15}$ FLOPs = $972$ petaFLOPs; it requires storing $\numdata$ vectors of size $\numparam \numtime$ which, assuming double precision, is $8$~TB of storage in total.
The one-time cost for computing the FFTs of the matrix blocks is negligible relative to the cost of pre-computing the operator.
The cost of performing FFT-accelerated matvecs of $\blocktoep$ and $\blocktoep^*$, which involves a number of different operations (see \Cref{sec:methods}), is approximately $8 \numparam \numdata \numtime = 8 \cdot 10^{12}$~FLOPs or $8$~teraFLOPs,\footnote{The costs are dominated by the SBGEMV (see \Cref{sec:results}, \cref{fig:SingleGPUScaling}). The coefficient $8$ comes from complex-valued operations (6 FLOPs per multiplication, 2 FLOPs per addition).}
which is $\sim$1000$\times$ more efficient than the conventional method for each $\blocktoep$ or $\blocktoep^*$ application.
Solving the inverse problem with the FFT-accelerated matvecs thus costs 972 petaFLOPs (one-time setup cost) plus $2r \cdot 8 \cdot 10^{12}$ = $1.6 \cdot 10^{18}$ FLOPs ($10^5$ FFT-based matvecs of $\blocktoep$ and $\blocktoep^*$), or 2.57 ExaFLOPs in total.\footnote{\rev{Each Hessian action requires a matvec with $\blocktoep$ and $\blocktoep^*$. For an iterative solver (e.g., CG) to converge on~\cref{eq:LinearInverseProblem}, approximately $r \sim 10^5$ iterations are required (using regularization preconditioning).}}
This is over $750 \times$ more efficient than solving the inverse
problem using the conventional way (i.e., based on forward/adjoint PDE
solves) of performing Hessian matvecs described
in \cref{sec:HessianMatvecforInverseProblems}.  
We emphasize that the realized efficiency gain will in many cases be much larger, particularly when (1) the PDE time-stepping method uses a much smaller time step for stability (CFL condition) or accuracy than the temporal discretization of parameters and data; 
(2) the PDE is spatially discretized with a high-order method; 
or (3) the LTI system describes a multiphysics or mixed problem with a larger number of variables (discretized with many DOFs per grid point).

\revstart
It is important to note that $\blocktoep$ is in general a dense matrix even though the spatially discretized PDE operator $\vb{A}$ may be sparse. While there may be some sparsity structure in $\blocktoep$ resulting from the specific location of the observation points in the domain and the speed of information propagation in the problem, this structure is highly problem-dependent and thus will not be considered in the development of the FFT-based matvec algorithm. The sparsity of $\vb{A}$ is exploited in the traditional PDE-based Hessian actions discussed previously. However, even though $\vb{A}$ is sparse, its size is $N_u \times N_u$, where $N_u$ is the spatial state dimension. In contrast, $\blocktoep$ entirely avoids the state dimension and instead deals with the spatial parameter dimension $N_m$, which is often much less than $N_u$. Even if $N_m \approx N_u$, methods that involve the PDE state directly will have to employ a time-stepping scheme that introduces further costs due to CFL constraints (explicit) or sparse linear solves (implicit). Moreover, while optimized sparse linear algebra packages exist for GPUs (e.g., cuSPARSE, cuDSS), sparse matrices have irregular memory access patterns that make it difficult to achieve high GPU performance. The FFT-based matvec algorithm relies entirely on dense linear algebra and FFT routines that map well onto GPUs. As a result, the performance of the FFT-based matvec is much greater than that of methods utilizing PDE solves to compute Hessian actions.  
\revend

The resulting savings from performing fast FFT-based Hessian matvecs are important for several reasons: (1) they can enable solving large-scale inverse problems that may otherwise be prohibitively expensive to solve; (2) they can significantly reduce the cost of inverting from many different data vectors, in which case the cost of pre-computing $\blocktoep$ and $\blocktoep^*$ is easily amortized; and (3) they can make real-time inversion feasible for time-sensitive applications.
Moreover, the FFT-based Hessian matvecs are well-suited for GPU-accelerated computation, even when the forward and adjoint problems employ implicit solvers, adaptive mesh refinement, low order discretization, or other methods that are typically not amenable to achieving peak performance on GPUs.

The remainder of this paper addresses the issue of efficiently applying $\blocktoep$ and $\blocktoep^*$ of the structure \eq{ShiftInvariance}.
Note that while the p2o map is shift-invariant in time (i.e.~$\blocktoep$ has Toeplitz structure with respect to time steps), the blocks $\blocktoeptxt_{ij}$ themselves are not assumed to have any special structure and, typically, they are not even square matrices ($\numparam \ne \numdata$).
We assume that \emph{either} $\blocktoep$ \emph{or} $\blocktoep^*$ is available in compactly stored form (see Section~\ref{sec:methods} for details).
For a discussion of the broader context of solving inverse problems
involving shift-invariant systems, we refer
to~\cite{henneking2025goal}. The present paper focuses on 
fast and scalable FFT-based algorithms for applying the p2o and
adjoint p2o operators to vectors, which are an essential building block for solving such inverse problems. Before describing our method in detail in \Cref{sec:methods}, we briefly introduce the well-known FFT-based matvec algorithm for Toeplitz matrices.

\subsection{FFT-Based Matvec for Toeplitz Matrices}\label{sec:FFTMatrixVectorProductforToeplitzMatrices}

This section describes the algorithm for computing matvecs of a generic Toeplitz matrix $\vb{M}_{\text{toep}} \in \R^{n \times n}$ given by
\begin{align}\label{eq:ToeplitzMatrix}
    \vb{M}_{\text{toep}} = \mqty[m_0&m_{-1}&\cdots&m_{-(n-1)}\\
    m_1&m_0&\ddots&\vdots\\
    \vdots&\vdots&\ddots&m_{-1}\\
    m_{n-1}&m_{n-2}&\cdots&m_0]_{\nonumber}
\end{align}
and a vector $\vb{x} \in \R^{n}$. \rev{Note that in this discussion the matrix $\vb{M}_{\text{toep}}$ is not related to a specific LTI system and its entries $m_{i}$ are not parameter values as in the previous discussions.}

To compute the matvec $\vb{M}_{\text{toep}}\vb{x}$, the matrix $\vb{M}_{\text{toep}}$ is first embedded inside the circulant matrix $\vb{M}_{\text{circ}} \in \R^{2n \times 2n}$ given by
\begin{equation}\label{eq:CirculantMatrix}
 \resizebox{\textwidth}{!}{$%
 \begin{aligned}[t]
	&\hspace{10pt} \vb{M}_{\text{circ}} := \mqty[\vb{M}_{\text{toep}}&\vb{M}'\\ \vb{M}'&\vb{M}_{\text{toep}}] := \\
	&\left[ \begin{array}{ccccc|ccccc}
    m_0&m_{-1}&\cdots&\cdots&m_{-(n-1)}&0&m_{n-1}&\cdots&m_2&m_1\\
                m_1&m_0&\ddots&&\vdots&m_{-(n-1)}&\ddots&\ddots&&m_2\\
                m_2&\ddots&\ddots&\ddots&\vdots&\vdots&\ddots&&\ddots&\vdots \\
                \vdots&&\ddots&m_0&m_{-1}&m_{-2}&&\ddots&\ddots&m_{n-1}\\[4pt]
                m_{n-1}&\cdots&\cdots&m_1&m_0&m_{-1}&m_{-2}&\cdots&m_{-(n-1)}&0\\[2pt]
                \hline
                0 &m_{n-1}&\cdots&\cdots&m_1&m_0&m_{-1}&\cdots&\cdots&m_{-(n-1)}\\
                m_{-(n-1)}&\ddots&\ddots&&m_2&m_1&m_0&\ddots&&\vdots\\
               \vdots&\ddots&&\ddots&\vdots&m_2&\ddots&\ddots&\ddots&\vdots \\
               m_{-2}&&\ddots&\ddots&m_{n-1}&\vdots&&\ddots&m_0&m_{-1}\\[4pt]
               m_{-1}&m_{-2}&\cdots&\cdots&0&m_{n-1}&\cdots&\cdots&m_1&m_0
     \end{array} \right] .
\end{aligned}$%
}
\end{equation}

It is well known~\cite{gray2006toeplitz} that the DFT matrix $\vb{D}$ diagonalizes the circulant matrix $\vb{M}_{\text{circ}}$; that is, $\vb{M}_{\text{circ}} = (2n)^{-1/2}\vb{D}^{-1}\text{diag}\qty(\qty(\widehat{\vb{M}}_{\text{circ}})_{0:})\vb{D}$, where $\qty(\vb{M}_{\text{circ}})_{0:}$ is the first column of $\vb{M}_{\text{circ}}$, and $\hat{}$ denotes the FFT.\@ So, to compute $\vb{M}_{\text{toep}}\vb{x}$, set $\vb{u} := \qty[\vb{x}\ \vb{0}]^T \in \R^{2n}$, and compute
\begin{align}\label{eq:FFTMatrixVectorProduct}
    \vb{M}_{\text{circ}}\vb{u} = \frac{1}{\sqrt{2n}}\vb{D}^{-1}\text{diag}\qty(\qty(\widehat{\vb{M}}_{\text{circ}})_{0:})\vb{D}\vb{u} = \text{IFFT}\qty(\frac{1}{\sqrt{2n}}\qty(\widehat{\vb{M}}_{\text{circ}})_{0:}\odot\hat{\vb{u}}),
\end{align}
where IFFT is the inverse FFT, and $\odot$ denotes the elementwise product. Observe that $\vb{M}_{\text{toep}}\vb{x}$ is the first block of the result.

%% file: methods.tex
This section describes the algorithm for computing the matvec with a block-triangular Toeplitz matrix $\blocktoep$. The algorithm for the case where $\blocktoep$ is block lower-triangular is considered first; the upper-triangular case, which can be handled similarly, is discussed later. Recalling the shift-invariant system \eq{ShiftInvariance}, the matrix $\blocktoep$ is structured as follows:
\begin{align}\label{eq:BlockToeplitzMatrix}
    \blocktoep =
    \left[ \begin{array}{ccccc}
	\blocktoeptxt_{11} & 0 & 0 & \cdots & 0 \\[2pt]
	\blocktoeptxt_{21} & \blocktoeptxt_{11} & 0 & \cdots & 0 \\
	\blocktoeptxt_{31} & \blocktoeptxt_{21} & \blocktoeptxt_{11} & \ddots & \vdots \\
	\vdots & \vdots & \ddots & \ddots & 0 \\[2pt]
	\blocktoeptxt_{\numtime,1} & \blocktoeptxt_{\numtime-1,1} & \cdots & \blocktoeptxt_{21} & \blocktoeptxt_{11} 
	\end{array} \right] .
\end{align}
$\blocktoep$ has block dimension $\numtime \times \numtime$ and $\blocktoeptxt_{ij}\in \R^{\numdata \times \numparam}$. Recall the definitions of the sizes:

\begin{itemize}
	\item $\numtime$ is the number of time steps;
	\item $\numdata$ is the spatial dimension of the data (e.g.~number of sensors);
	\item $\numparam$ is the spatial dimension of the parameters (e.g.~number of sources).
\end{itemize}

The algorithm for computing the matvec with the matrix $\blocktoep$ can be broken down into two steps: 1) a setup phase where the matrix $\blocktoep$ is read from file and transformed into a new matrix $\toepblock$ with triangular Toeplitz blocks, and 2) a matvec phase where the matvec is computed using the FFT.\@ The setup and $\blocktoep$ matvec algorithms are given in~\cref{alg:setup,alg:matvec}, respectively. In~\cref{sec:MatvecWithBlockToeplitzConjugateTranspose}, the algorithm for computing the matvec with the matrix $\blocktoep^*$ (which can also be applied to general block upper-triangular Toeplitz matrices) is described. The algorithm for computing the $\blocktoep^*$ matvec is given in~\cref{alg:matvec2}. Throughout these algorithms, several index transformations are required; these, along with the notation used in the remainder of the paper, will be discussed first.

\subsection{Index Transformations and Notation}\label{sec:IndexTransformation}
The matrix $\blocktoep$ is block-triangular Toeplitz when written with \textit{time-outer-space-inner} (TOSI) ordering. That is, blocks of $\blocktoep$ correspond to time steps, and each block of $\blocktoep$ corresponds to spatial information of the data and parameters. This index ordering can be switched to \textit{space-outer-time-inner} (SOTI) ordering. In SOTI form, $\blocktoep$ has $\numdata \times \numparam$ blocks; \textit{each block} is lower-triangular Toeplitz and has size $\numtime \times \numtime$. The TOSI and SOTI orderings also apply to the discretized parameter vector $\paramvec$ and data vector $\datavec$. Different steps of the matvec algorithms involve quantities represented in either form. To denote quantities in the SOTI ordering, we use a tilde; quantities in TOSI ordering are denoted without tildes. For example, the discrete p2o map is represented by $\blocktoep$ in the TOSI ordering and $\toepblock$ in the SOTI ordering. They are related by $(\toepblocktxt_{ij})_{kl} = (\blocktoeptxt_{kl})_{ij}$. Similarly, $(\tparamvec_j)_l = (\paramvec_l)_j$ and $(\tdatavec_i)_k = (\datavec_k)_i$. The equation~\cref{eq:p2o} is invariant under this change of index: $\toepblock \tparamvec = \tdatavec$. In the implementation, where the matrices and vectors are stored as 1-dimensional arrays, changing between TOSI and SOTI ordering (for local quantities) corresponds to a transpose or swapaxes operation. GPU algorithms for these operations have been well studied~\cite{jodra2015efficient,ruetsch2009optimizing}, and we can use them directly (e.g.\~from the cuTENSOR\footnote{\url{https://docs.nvidia.com/cuda/cutensor/latest/index.html}} library). Note that in an expression such as $(\blocktoeptxt_{kl})_{ij}$, the indices inside the parentheses denote ``outer'' indices (i.e., block indices $(k,l)$), \revstart and the indices outside the parentheses denote ``inner'' indices (i.e., element indices $(i,j)$ within block $(k,l)$).

These notational conventions are summarized in~\cref{tab:notation}.\revend

The matvec algorithm is designed with a multi-GPU implementation in mind. As such, the matrices and vectors are partitioned over the processors. In most cases, whether an expression refers to a local or global quantity should be evident from context. In places where there may be ambiguity, we use the notation $\vb{v}^{\text{G}}$ to denote a global quantity and $\vb{v}^{\text{L};ij}$ to denote a local quantity on the processor $\proc_{ij}$. 

The Fourier transform is used throughout the algorithm. Fourier transformed quantities are denoted with hats (e.g. $\widehat{\vb{v}}$). 

\revstart

\begin{table}[htbp]
\revstart
\centering
\resizebox{\textwidth}{!}{%
\begin{tabular}{@{}lll@{}}
\toprule
Notation          & Meaning                         & Example                                  \\ \midrule
Standard              & TOSI (time-outer-space-inner indexing) & $\blocktoep, \paramvec, \datavec$                              \\
Tilde              & SOTI (space-outer-time-inner indexing) & $\toepblock, \tparamvec, \tdatavec$                             \\
Hat &
  Fourier space matrix/vector &
  $\widehat{\blocktoep}$, $\widehat{\paramvec}$, $\widehat{\datavec}$ \\
Superscript $\text{L}; ij$ &
  Local portion of matrix/vector on processor $p_{ij}$ &
  $\blocktoep^{\text{L};ij}$, $\paramvec^{\text{L};ij}$ \\
Superscript $\text{G}$ &
  Global matrix/vector (omitted if clear from context) &
  $\datavec^{\text{G}}$, $\paramvec^{\text{G}}$ \\
$(\square_{ij})_{kl}$ &
  Block $(i,j)$ element $(k,l)$ of $\square$ &
  $\qty(\blocktoep_{ij})_{kl}$, $\qty(\paramvec_{i})_{k}$, $({\datavec}_{j})_l$ \\ \bottomrule
\end{tabular}%
}
\caption{\rev{Notational conventions. Stacked notations can result in expressions like $\widehat{\toepblock}{}^{\text{L};ij}$\hspace{-3pt}.}}
\label{tab:notation}
\revend
\end{table}
\revend

\subsection{Preprocessing and Partitioning of the Matrix}\label{sec:ReadingBlockToeplitz}
As mentioned in~\cref{sec:InverseProblemsInvolvingShiftInvariantSystems}, the matrix $\blocktoep$ is formed by computing $\numparam$ forward solves or $\numdata$ adjoint solves. When performing these solves, it is natural to work in TOSI ordering. Thus, before computing matvecs, the matrix data has to be preprocessed. The first preprocessing step is to convert the matrix to SOTI ordering. This is most easily done via a short Python script (using $\texttt{numpy.reshape}$, for example). 

Once the matrix data is in SOTI ordering, it can be partitioned among the processors (GPUs) in the multi-GPU setup. For the distributed-memory parallelism model, a 2D grid of processors of size $\numprocrows \times \numproccols$ is assumed. In SOTI ordering, $\toepblock$ is partitioned along the outer indices; each processor $\proc_{ij}$, $0\leq i \leq \numprocrows-1$, $0\leq j\leq \numproccols-1$, holds the $\procrowblock\times\proccolblock$ sub block matrix shown in~\cref{eq:ProcToeplitzBlocks} (disregarding edge cases). Here, $\procrowblock = \lceil{\frac{\numdata}{\numprocrows}}\rceil$ and $\proccolblock = \lceil{\frac{\numparam}{\numproccols}}\rceil$. For ease of presentation, we will only show examples of the partitioning where $\numdata/\numprocrows$ and $\numparam/\numproccols$ are whole numbers.

\begin{align}\label{eq:ProcToeplitzBlocks}
    \toepblock^{\text{L};ij} = \mqty[
        \toepblocktxt_{i\procrowblock,j\proccolblock} & \toepblocktxt_{i\procrowblock,j\proccolblock+1} & \cdots & \toepblocktxt_{i\procrowblock,(j+1)\proccolblock-1}\\
        \toepblocktxt_{i\procrowblock+1,j\proccolblock} & \toepblocktxt_{i\procrowblock+1,j\proccolblock+1} & \cdots & \toepblocktxt_{i\procrowblock+1,(j+1)\proccolblock-1}\\
        \vdots & \vdots & \ddots & \vdots\\
        \toepblocktxt_{(i+1)\procrowblock-1,j\proccolblock} & \toepblocktxt_{(i+1)\procrowblock-1,j\proccolblock+1} & \cdots & \toepblocktxt_{(i+1)\procrowblock-1,(j+1)\proccolblock-1}\\
    ].
\end{align}

When computing matvecs, the parameter and data vectors should also be in SOTI ordering. In this case, the parameter vector is partitioned along the first row of processors $\proc_{0,j}$, and the data vector is partitioned along the first column of processors $\proc_{i,0}$.~\Cref{eq:ProcParamVec,eq:ProcDataVec} show an example of the vector partitioning:
\begin{align}
 \tparamvec^{\text{L};0,j} &=\mqty[
    \tparamtxt_{j\proccolblock}&
    \tparamtxt_{j\proccolblock + 1}&
    \hdots&
    \tparamtxt_{(j+1)\proccolblock - 1}
]^T,\label{eq:ProcParamVec}\\
 \tdatavec^{\text{L};i,0} &= \mqty[
    \tdatatxt_{i\procrowblock}&
    \tdatatxt_{i\procrowblock + 1}&
    \hdots&
    \tdatatxt_{(i+1)\procrowblock - 1}
]^T.\label{eq:ProcDataVec}
\end{align}

The index-transformed matvec $\tdatavec = \toepblock\tparamvec$ can be computed as in~\cref{eq:BlockMatvec}:
\begin{align}\label{eq:BlockMatvec}
    \tdatavec_k^{\text{G}} = \qty(\toepblock^{\text{G}}\tparamvec^{\text{G}})_k &= \sum_{k=0}^{\numparam-1} \toepblock_{kl}^{\text{G}} \tparamvec_l^{\text{G}}\\
    \Rightarrow \tdatavec_k^{\text{L};i,0} &=\sum_{j=0}^\numproccols\sum_{k=0}^{\proccolblock-1} \toepblock_{kl}^{\text{L};ij} \tparamvec_l^{\text{L};ij}.
\end{align}

First, the parameter vector is broadcast down each processor column so that $\tparamvec^{\text{L};ij} = \tparamvec^{\text{L};0,j}$ for all $0\leq i \leq \numprocrows$.
Then, each processor $\proc_{ij}$ computes the local matvecs $\toepblock_{kl}^{\text{L};ij} \tparamvec_l^{\text{L};ij}$. Finally, a reduction with summation is computed over each row of the processor grid.

The algorithm for matvecs with $\blocktoep^*$ is structured in the same way, as shown in~\cref{eq:BlockTransMatvec}. First, the data vector is broadcast down each processor row so that $\tdatavec^{\text{L};ij} = \tdatavec^{\text{L};i,0}$ for all $0\leq j \leq \numprocrows$. Then, each processor $\proc_{ij}$ computes the local matvecs $\toepblock_{kl}^{*;\text{L};ij} \tdatavec_l^{\text{L};ij}$. Finally, a reduction with summation is computed over each column of the processor grid.

\begin{align}\label{eq:BlockTransMatvec}
    \tparamvec_k^{\text{G}} = \qty(\toepblock^{*;\text{G}}\tdatavec^{\text{G}})_k &= \sum_{k=0}^{\numdata-1} \toepblock_{kl}^{*;\text{G}} \tdatavec_l^{\text{G}}\\
    \Rightarrow \tparamvec_k^{\text{L};0,j} &=\sum_{i=0}^\numprocrows\sum_{k=0}^{\procrowblock-1} \toepblock_{kl}^{*;\text{L};ij} \tdatavec_l^{\text{L};ij}.
\end{align}

In either case, each local block matvec (e.g.~$\toepblock_{kl}^{\text{L};ij} \tparamvec_l^{\text{L};ij}$) involves a lower-triangular Toeplitz matrix. The block lower-triangular Toeplitz structure of $\blocktoep$ (in the TOSI ordering) translates to the lower-triangular Toeplitz structure of each block of $\toepblock$ (in the SOTI ordering). Recall that we only compute/store the first block column of $\blocktoep$. This corresponds to only storing the first column of each block of $\toepblock$.

The next preprocessing step for $\toepblock$ is to pad the first column of each of its blocks and take the Fourier transform (as explained in~\cref{sec:FFTMatrixVectorProductforToeplitzMatrices}). This can be done via batched FFT methods on GPUs (we use the cuFFT library\footnote{\url{https://docs.nvidia.com/cuda/cufft/index.html}}). We denote by $\widehat{\toepblock}$ the result of the aforementioned process. The final preprocessing step is to convert $\widehat{\toepblock}$ back to TOSI ordering; the result is denoted by $\widehat{\blocktoep}$. In SOTI ordering, $\widehat{\toepblock}$ consists of $\numdata \times \numparam$ diagonal blocks. A possible way to compute a matvec from there would be to take elementwise products and then sum over the results of each block row. This is the most natural way to formulate the algorithm, as it parallels the process in~\cref{sec:FFTMatrixVectorProductforToeplitzMatrices}. However, $\widehat{\blocktoep}$ is a block-diagonal matrix; after performing the corresponding reordering to the vector, the local matvec then becomes a matvec between a block diagonal matrix and vector (both consisting of complex numbers). This algorithm has the advantage that all arithmetic operations involve data that is contiguous in memory. Moreover, it can leverage routines from libraries such as cuBLAS that are well known to achieve high performance on GPUs~\cite{relton2016comparison}. In contrast, the former algorithm involves strided memory access patterns and has to be implemented using custom GPU kernels (see~\Cref{sec:appendix-alt-matvec} for a discussion of this algorithm).

The conversion of local portions of $\widehat{\toepblock}$ from SOTI ordering back to TOSI ordering can be achieved with a swapaxes kernel~\cite{jodra2015efficient}.\footnote{To be used in cuBLAS operations, the matrix should be stored in column-major order. Initially, the SOTI matrices are stored in row-major order. However, the swapaxes operation accomplishes both the SOTI to TOSI conversion and the switch to column-major ordering in the same kernel. Column-major storage is not needed anywhere else in the algorithm, so all other arrays are assumed to be row-major. \rev{This operation is inexpensive compared to the cost of computing the matrix FFT.}} The entries of the result, $\widehat{\blocktoep}$, can then be stored for later use. If the same processor partitioning is used for this preprocessing/setup phase and the later matvec computations, no further reordering operations need to be done. Otherwise, the entries of $\widehat{\blocktoep}$ need to be rearranged to match whatever partitioning of $\toepblock$ is used in the matvec computations. Again, this is most easily done with a simple Python script. Note that all reordering operations (e.g.~\texttt{TOSI\_TO\_SOTI} and \texttt{SOTI\_TO\_TOSI}) are only performed on \textit{local} portions of matrices and vectors; no distributed communication is required.

Throughout the rest of this paper, we assume that the same partitioning is used for the preprocessing/setup phases and the matvec computations. \Cref{alg:setup,alg:matvec,alg:matvec2} detail the matvec algorithms. A discussion of each step of the algorithms can be found in subsequent sections.

\begin{algorithm}
\caption{Setup phase for computing the matvec with a block-triangular Toeplitz matrix $\blocktoep$ and $\blocktoep^*$.}\label{alg:setup}
\begin{algorithmic}[1]
\Procedure{setup}{\texttt{filename}, $\numparam$, $\numtime$, $\numdata$, $\numprocrows$, $\numproccols$} \Comment{2D processor index $(i,j)$}
\State{$\toepblock \gets$ \texttt{TOSI\_to\_SOTI}$\qty(\blocktoep)$} \Comment{~\cref{sec:ReadingBlockToeplitz,sec:IndexTransformation}; done separately prior to reading from file.}
\State{$\procrowblock \gets \lceil{\frac{\numdata}{\numprocrows}}\rceil$; $\proccolblock \gets \lceil{\frac{\numparam}{\numproccols}}\rceil$}
\State{$\toepblock_{kl}^{\text{L};ij} \gets$ \texttt{read}(\texttt{filename}, \texttt{row\_start=} $i \procrowblock$, \texttt{col\_start=} $j \proccolblock$}\Comment{~\cref{sec:ReadingBlockToeplitz}; assumes SOTI ordering}
\State{$\widehat{\toepblock}{}^{\text{L};ij}\gets$ \texttt{batched\_padded\_fft}$\qty(\toepblock_{kl}^{\text{L};ij})$} \Comment{~\cref{sec:LocalMatvec}\\}
\State{$\widehat{\blocktoep}^{\text{L};ij}\gets$ \texttt{SOTI\_to\_TOSI}$\qty(\widehat{\toepblock}{}^{\text{L};ij})$} \Comment{~\cref{sec:IndexTransformation}; add zero padding}

\noindent\Return{$\widehat{\blocktoep}$}
\Comment{~\cref{sec:ReadingBlockToeplitz,sec:IndexTransformation}}
\EndProcedure{}
\end{algorithmic}
\end{algorithm}

\begin{algorithm}
\caption{Matvec phase for computing the matvec with the block lower-triangular Toeplitz matrix $\blocktoep$. Assumes $\tparamvec$ in SOTI ordering (else run \texttt{TOSI\_to\_SOTI} first). Returns $\tdatavec$ in SOTI ordering. $\widehat{\blocktoep}$ comes from running \texttt{setup}.}\label{alg:matvec}
\begin{algorithmic}[1]
\Procedure{matvec}{$\widehat{\blocktoep}$, $\tparamvec$, $\numparam$, $\numtime$, $\numdata$, $\numprocrows$, $\numproccols$}\Comment{2D processor index $(i,j)$}
\State{$\tparamvec_l^{\text{L};ij} \gets$ \texttt{broadcast}$\qty(\tparamvec)$} \Comment{broadcast input vector to all processor rows}
\State{$\widehat{\tparamvec}{}^{\text{L};ij}_l\gets$ \texttt{batched\_padded\_fft}$\qty(\tparamvec_l^{\text{L};ij})$} \Comment{~\cref{sec:LocalMatvec}}
\State{$\widehat{\paramvec}^{\text{L};ij}\gets$ \texttt{SOTI\_to\_TOSI}$\qty(\widehat{\tparamvec}{}^{\text{L};ij})$} \Comment{~\cref{sec:LocalMatvec}}
\State{$\widehat{\datavec}^{\text{L};ij} \gets$ \texttt{apply\_matrix}$\qty(\widehat{\blocktoep}^{\text{L};ij}, \widehat{\paramvec}^{\text{L};ij})$} \Comment{~\cref{sec:LocalMatvec}}
\State{$\widehat{\tdatavec}{}^{\text{L};ij}\gets$ \texttt{TOSI\_to\_SOTI} $\qty(\widehat{\datavec}^{\text{L};ij})$} \Comment{~\cref{sec:LocalMatvec}}
\State{$\tdatavec^{\text{L};ij}\gets$ \texttt{batched\_ifft\_and\_unpad}$\qty(\widehat{\tdatavec}{}^{\text{L};ij})$}  \Comment{~\cref{sec:LocalMatvec}}
\State{$\tdatavec^{\text{G}}\gets$ \texttt{reduction}$\qty(\tdatavec^{\text{L};ij})$}\Comment{~\cref{sec:ReductionOverMultipleGPUs}; reduce each processor row}

\noindent\Return{$\tdatavec^{\text{G}}$} \Comment{SOTI ordering ~\cref{sec:IndexTransformation}}
\EndProcedure{}
\end{algorithmic}
\end{algorithm}

\begin{algorithm}
\caption{Matvec phase for computing the matvec with the block lower-triangular Toeplitz matrix $\blocktoep^*$. Assumes $\tdatavec$ in SOTI ordering (else run \texttt{TOSI\_to\_SOTI} first). Returns $\tparamvec$ in SOTI ordering. $\widehat{\blocktoep}$ comes from running \texttt{setup}. No need to rerun setup if already run for $\blocktoep$ matvecs.}\label{alg:matvec2}
\begin{algorithmic}[1]
\Procedure{matvec}{$\widehat{\blocktoep}$, $\tparamvec$, $\numparam$, $\numtime$, $\numdata$, $\numprocrows$, $\numproccols$}\Comment{2D processor index $(i,j)$}
\State{$\tdatavec_l^{\text{L};ij} \gets$ \texttt{broadcast}$\qty(\tdatavec)$} \Comment{broadcast input vector to all processor columns}
\State{$\widehat{\tdatavec}{}^{\text{L};ij}_l\gets$ \texttt{batched\_padded\_fft}$\qty(\tdatavec_l^{\text{L};ij})$} \Comment{~\cref{sec:LocalMatvec}}
\State{$\widehat{\datavec}^{\text{L};ij}\gets$ \texttt{SOTI\_to\_TOSI}$\qty(\widehat{\tdatavec}{}^{\text{L};ij})$} \Comment{~\cref{sec:LocalMatvec}}
\State{$\widehat{\paramvec}^{\text{L};ij} \gets$ \texttt{apply\_matrix}$\qty(\widehat{\blocktoep}^{\text{L};ij}, \widehat{\datavec}^{\text{L};ij},\texttt{conjugate=TRUE})$} \Comment{~\cref{sec:LocalMatvec}}
\State{$\widehat{\tparamvec}{}^{\text{L};ij}\gets$ \texttt{TOSI\_to\_SOTI} $\qty(\widehat{\paramvec}^{\text{L};ij})$} \Comment{~\cref{sec:LocalMatvec}}
\State{$\tparamvec^{\text{L};ij}\gets$ \texttt{batched\_ifft\_and\_unpad}$\qty(\widehat{\tparamvec}{}^{\text{L};ij})$}  \Comment{~\cref{sec:LocalMatvec}}
\State{$\tparamvec^{\text{G}}\gets$ \texttt{reduction}$\qty(\tparamvec^{\text{L};ij})$}\Comment{~\cref{sec:ReductionOverMultipleGPUs}; reduce each processor row}

\noindent\Return{$\tparamvec^{\text{G}}$} \Comment{SOTI ordering ~\cref{sec:IndexTransformation}}
\EndProcedure{}
\end{algorithmic}
\end{algorithm}

\subsection{Local Matvecs}\label{sec:LocalMatvec}

After partitioning the matrix $\blocktoep$ as discussed in~\cref{sec:ReadingBlockToeplitz}, the problem has been reduced to computing the local matvecs $\toepblock_{kl}^{\text{L};ij} \tparamvec_l^{\text{L};ij}$ (or correspondingly $\toepblock_{kl}^{*;\text{L};ij} \tdatavec_l^{\text{L};ij}$). As discussed in~\cref{sec:ReadingBlockToeplitz}, when the matrices and vectors are Fourier transformed and converted back to TOSI ordering, this operation is equivalent to applying a block diagonal matrix to a vector. 

Given a vector $\tparamvec$ or $\tdatavec$ in SOTI ordering, The first step is to pad each block with zeros and take the Fourier transform. Padding is computed with trivial custom CUDA kernels, and the batched FFTs are computed with cuFFT. After this step, the vectors are reordered to TOSI ordering.

The next step is to apply the block diagonal local matrix $\widehat{\blocktoep}^{\text{L}}$ to the local vector $\widehat{\paramvec}^{\text{L}}$ (or correspondingly, $\widehat{\blocktoep}^{*;\text{L}}$ to the vector $\widehat{\datavec}^{\text{L}}$). This is achieved via the cuBLAS\footnote{\url{https://docs.nvidia.com/cuda/cublas/}} \texttt{gemvStridedBatched} operation. For the case of $\widehat{\blocktoep}^{*;\text{L}}$ applications, the conjugate transpose of each of the diagonal blocks of $\widehat{\blocktoep}^{\text{L}}$ is applied instead. This transpose operation is not explicit; it is implicitly applied by the cuBLAS kernel when given the appropriate parameter. See~\cref{sec:MatvecWithBlockToeplitzConjugateTranspose} for a discussion on why only the matrix $\blocktoep$ needs to be partitioned and stored. 

After the matrix application, the resulting vector is converted back to SOTI ordering. Then, a batched IFFT is applied to this vector, and the result is unpadded (again with a trivial custom CUDA kernel). The result of this unpadding operation is the output: $\tdatavec^{\text{L}}$ for the $\blocktoep$ matvec and $\tparamvec^{\text{L}}$ for the $\blocktoep^*$ matvec. Again, note that all reordering operations are only applied to \textit{local} portions of the vectors, so no distributed communication is necessary.

\subsection{Matvecs With the Transposed Matrix}\label{sec:MatvecWithBlockToeplitzConjugateTranspose}

Next, matvecs with the matrix $\blocktoep^*$, which is a block \textit{upper}-triangular Toeplitz matrix, are described. One way to proceed would be to use \rev{\cref{alg:matvec}} as before, but modify the padding before taking the FFT.\@ However, this would require storing local parts of $\blocktoep$ and $\blocktoep^*$ separately, and may also require repartitioning. To avoid this extra cost, recall a key property of the Fourier transform that allows us to use the same algorithm as before with very little modification.

Recall from~\cref{sec:FFTMatrixVectorProductforToeplitzMatrices} that the first step to computing the matvec of a Toeplitz matrix $\vb{M}_{\text{toep}}$ with a vector $\vb{x}$ is to form the circulant matrix $\vb{M}_{\text{circ}}$ (cf.~\cref{eq:CirculantMatrix}). Furthermore, recall the relation $\vb{M}_{\text{circ}} = (2n)^{-1/2}\vb{D}^{-1}\text{diag}\qty(\qty(\widehat{\vb{M}}_{\text{circ}})_{0:})\vb{D}$, where $\vb{D}$ is the DFT matrix and $\square_{0:}$ refers to the first column of the matrix. However, note that since the Fourier transform is a unitary operation, $\vb{D}^{-1} = \vb{D}^*$. This implies 
\begin{equation}\label{eq:AdjointFFTMatvec}
\begin{split}
	\vb{M}_{\text{circ}}^* = \mqty[\vb{M}_{\text{toep}}^*&\qty(\vb{M}')^*\\
	\qty(\vb{M}')^*&\vb{M}_{\text{toep}}^*] = \frac{1}{\sqrt{2n}}\vb{D}^{-1}\text{diag}\qty(\qty(\widehat{\vb{M}}_{\text{circ}})_{0:})^*\vb{D} .
\end{split}
\end{equation}
Therefore, to compute the matvec with $\blocktoep^*$, the only required modification of~\rev{\cref{alg:matvec}} is taking the complex conjugate transpose of the blocks of $\widehat{\blocktoep}^{\text{L}}$ before applying it to $\widehat{\datavec}$. As mentioned in~\cref{sec:LocalMatvec}, the conjugate transpose operation is implicitly applied by cuBLAS. This is the only change required to compute the matvec with $\blocktoep^*$, and no extra storage or partitioning is necessary. Any generic block upper-triangular Toeplitz matrix can be handled in the same way.

\subsection{Reduction Over Multiple GPUs}\label{sec:ReductionOverMultipleGPUs}
After each processor (GPU) has computed its local matvec, the next step is to compute a reduction over either the rows ($\blocktoep$ matvec) or columns ($\blocktoep^*$ matvec) of the 2D processor grid. For communication between processors, NVIDIA's \texttt{NCCL} library is used to pass messages directly between GPUs, avoiding the need for copying data to the CPU and back. Specifically, the \texttt{ncclReduce} operation can be used for this task.

Note that to achieve the best performance over multiple consecutive matvecs, all memory and other structures needed for the computation (e.g.~temporary arrays, FFT plans) are allocated during the setup phase and reused over the multiple matvecs. After all matvecs are complete, this memory is freed. CUDA streams are also used to reduce kernel launch latency.

\subsection{Computational Cost}\label{sec:ComputationalCost}
Next, the theoretical time complexities for each of the steps in computing $\blocktoep$ and $\blocktoep^*$ matvecs listed in~\cref{sec:ReductionOverMultipleGPUs} are derived. The steps of the algorithm and their computational complexities (per GPU) are listed in~\cref{tab:ComputationalCost}. The complexities of reading the matrix from file, applying the index transformation, and moving data to the GPU are not calculated since these are one-time setup costs. The complexity of the FFT for the matrix is reported, but it is also a one-time setup cost. The complexity of each computational step is easily derived. The dimensions $\numparam$ and $\numdata$ are split among the $\numprocrows$ rows and $\numproccols$ columns of processors, respectively, giving the local sizes $\procrowblock$ and $\proccolblock$. The non-parallelizable dimension is $\numtime$. The FFTs scale as $2\numtime \log\qty(2\numtime)$ (factor 2 due to padding), and everything else scales as $\numtime$. There is a factor of 8 in the communication steps since double-precision floating point numbers are used.

\begin{table}[htbp]
    \centering
    \resizebox{\textwidth}{!}{%
    \begin{tabular}{@{}p{.255in}p{1.5in}p{1.5in}p{1.7in}@{}}
    \toprule
      & Description & Computational Cost (Total work per GPU) & Notes \\ \midrule
     \multirow{2}{*}[-0.5ex]{\rotatebox[origin=c]{90}{Setup}}& FFT Matrix & 
     $\order{2\procrowblock\proccolblock\numtime\log\qty(2\numtime)}$ & One time setup cost \\
     & \texttt{SOTI\_TO\_TOSI} & --- & One time memory operation\\ \\[-3pt]
     \hline
     \multirow{9}{*}[-0.5ex]{\rotatebox[origin=c]{90}{$\blocktoep$ Matvec}}& Broadcast Vector & $\order{\qty(\ell + 8\numtime\proccolblock/\beta)\log\numproccols}$ & Latency $\ell$; Bandwidth $\beta$\\
     & Pad Vector & $\order{2\proccolblock\numtime}$ & --- \\
     & FFT Vector & $\order{2\proccolblock\numtime\log\qty(2\numtime)}$ & --- \\
     & \texttt{SOTI\_TO\_TOSI} & --- & Memory operation \\
     & Apply Matrix & $\order{\procrowblock\proccolblock(\numtime + 1)}$ & Complex arithmetic\\
     & \texttt{TOSI\_TO\_SOTI} & --- & Memory operation\\
     & IFFT & $\order{2\procrowblock\numtime\log\qty(2\numtime)}$  & --- \\
     & Unpad Vector & $\order{2\procrowblock\numtime}$ & --- \\
     & Reduce Vector & $\order{\qty(\ell + 8\numtime\procrowblock/\beta)\log\numproccols}$ & Latency $\ell$; Bandwidth $\beta$ \\\hline
    \multirow{9}{*}[-0.5ex]{\rotatebox[origin=c]{90}{$\blocktoep^*$ Matvec}} & Broadcast Vector & $\order{\qty(\ell + 8\numtime\procrowblock/\beta)\log\numproccols}$ & Latency $\ell$; Bandwidth $\beta$ \\
    & Pad Vector & $\order{2\procrowblock\numtime}$ & --- \\
     & FFT Vector & $\order{2\procrowblock\numtime\log\qty(2\numtime)}$ & --- \\
    & \texttt{SOTI\_TO\_TOSI} & --- & Memory operation \\
     & Apply Matrix & $\order{\procrowblock\proccolblock(\numtime + 1)}$ & Complex arithmetic \\
     & \texttt{TOSI\_TO\_SOTI} & --- & Memory operation\\
     & IFFT & $\order{2\proccolblock\numtime\log\qty(2\numtime)}$ & --- \\
     & Unpad Vector & $\order{2\proccolblock\numtime}$ & --- \\
     & Reduce Vector & $\order{\qty(\ell + 8\numtime\proccolblock/\beta)\log\numprocrows}$ & Latency $\ell$; Bandwidth $\beta$ \\ \bottomrule
    \end{tabular}%
    }
    \caption{Steps of matvec algorithm for $\blocktoep$ and $\blocktoep^*$ matvecs along with their computational complexities per GPU.}\label{tab:ComputationalCost}
    \end{table}

The computational costs given in~\cref{tab:ComputationalCost} show that as the local sizes $\procrowblock$ and $\proccolblock$ increase (with $\numtime$ fixed, the matrix application \rev{via a Strided-Batched-GEneral-Matrix-Vector-product} (SBGEMV) cuBLAS kernel is asymptotically the most expensive step. 

\revstart
The dimension along which~\cref{alg:matvec} cannot be fully parallelized is $\numtime$. This algorithm scales as $\order{2\numtime\log\qty(2\numtime)}$ in this dimension. In contrast, a direct dense matvec would scale as $\order{\numtime^2/2}$ in the time dimension (factor 1/2 due to the block-triangular structure).
\revend

\subsection{{Communication-Aware Partitioning}}\label{sec:gpu-grid}
{Communication costs during the broadcast and reduction phases of the $\blocktoep$~and $\blocktoep^*$ matvecs present the primary barrier to the weak scalability of the matvec algorithms. To overcome the above we propose a communication-aware partitioning scheme. We use the asymptotic total communication cost $C(r, c)$ as an approximate cost model in our partitioning scheme.}
From~\cref{tab:ComputationalCost}, the total communication cost $C(r, c)$ for one matvec with $\blocktoep$ and one with $\blocktoep^*$ {is given below}. 
\begin{align}
     C(\numprocrows, \numproccols) &:= \qty(\ell + \frac{8\numtime\proccolblock}{\beta})\log\numprocrows + \qty(\ell + \frac{8\numtime\procrowblock}{\beta})\log\numproccols\\
     &= \qty(\ell + \frac{8\numtime\numparam}{\beta\numproccols})\log\numprocrows + \qty(\ell + \frac{8\numtime\numdata}{\beta\numprocrows})\log\numproccols
\end{align}
Thus, { for a given partition $p=\numprocrows\times \numproccols$, we can choose $\numprocrows$ and $\numproccols$ that minimize $C(\numprocrows, \numproccols)$.} We can rewrite $C(\numprocrows, \numproccols)$ as 
\begin{align}
    C(\numprocrows, \numproccols) = \ell \log p + \frac{8\numtime}{\beta}\qty(\frac{\numparam}{\numproccols}\log \numprocrows + \frac{\numdata}{\numprocrows}\log \numproccols).
\end{align}
In {the above} expression, the latency term $\ell \log p$ and the bandwidth term $8\numtime/\beta$ are just constant shifts or scaling of the total communication cost. As a result, it is sufficient to minimize the modified asymptotic cost function $\tilde{C}(\numprocrows, \numproccols)$ defined as 
\begin{align}
        \tilde{C}(\numprocrows, \numproccols) := \frac{\numparam}{\numproccols}\log \numprocrows + \frac{\numdata}{\numprocrows}\log \numproccols.
\end{align}
Rewriting {the above} expression using $p = \numprocrows \times \numproccols$ gives
\begin{align}
    \tilde{C}(r) &= \frac{\numparam \numprocrows}{p}\log \numprocrows + \frac{\numdata}{\numprocrows}\log\frac{p}{\numprocrows} = \numparam \qty(\frac{\numprocrows}{p}\log \numprocrows + \frac{10^l}{\numprocrows}\log\frac{p}{\numprocrows}),\label{eq:modcost}
\end{align}
where $l = \log_{10} \qty(\numdata / \numparam)$. Removing the overall scale factor $\numparam$, the modified asymptotic communication cost function $\tilde{C}$ is viewed as a function of the number of rows in the {partition} grid. The total number of processes $p$ as well as the log ratio $l$ of the global data dimension to global parameter dimension act as hyperparameters for $\tilde{C}$, now defined as
\begin{align}
    \tilde{C}(r;p,l) &:= \frac{\numprocrows}{p}\log \numprocrows + \frac{10^l}{\numprocrows}\log\frac{p}{\numprocrows}.\label{eq:comm-cost}
\end{align}
For a given value of $p$ and $l$ (generally known before any computations), the theoretical optimal grid configuration is the one with a number of rows that minimizes the modified asymptotic communication cost in~\cref{eq:comm-cost}. Now, there are some extra considerations: the number of processor rows $\numprocrows$ and columns $\numproccols$ both have to be \textit{natural numbers} multiplying to $p$. So, in practice, one can minimize the real-valued function in~\cref{eq:comm-cost} and then choose the best value for $\numprocrows$ near the real number minimizer that leads to a valid processor grid. 

Another consideration that is important in the practical implementation is that of \textit{on-node} vs. \textit{off-node} communication costs. Communication costs between processes on the same node are generally much smaller than between processes on different nodes. For a cluster where there are $k$ GPUs per node, communication costs along columns of the processor grid likely remain unchanged if rows are added so that $k | \numprocrows$. The exception to this is if $\numprocrows$ is changed from 1 to $k$, in which case communication costs will definitely increase along columns. As a result, the procedure in~\cref{alg:grid} is recommended to determine the optimal value of $\numprocrows$. It is important to verify the optimality of the configuration through testing (at least by slightly perturbing the value of $\numprocrows$) as hardware and network configurations can have a substantial impact on the optimal processor grid shape.

\begin{algorithm}
\caption{Selecting Optimal Processor Grid Shape\\ ($p = $ total number of processors; $l = \log_{10} \qty(\numdata / \numparam)$; $k$ GPUs per node)} \label{alg:grid}
\begin{algorithmic}[1]
\Procedure{select\_grid}{$p$, $l$, $k$}
\State{$\tilde{\numprocrows} \gets \argmin_\numprocrows \tilde{C}\qty(\numprocrows; p, l)$}
\If{$\tilde{\numprocrows} = 1$}
\Return{$(1, p)$}
\EndIf
\If{$\tilde{\numprocrows} = p$}
\Return{$(p, 1)$}
\EndIf
\State{Choose $\numprocrows\in \N$ close to $\tilde{\numprocrows}$ and $\numproccols \in \N$ close to $p/\tilde{\numprocrows}$ so that $p = \numprocrows\numproccols$. Try to ensure $k | \numprocrows$ and/or $k |\numproccols$. If $l \geq 0$, choose $\numprocrows \geq \numproccols$. Else choose $\numprocrows < \numproccols$.}
\Statex{}
\Return{$\numprocrows, \numproccols$} \Comment{Check optimality of grid in practice.}
\EndProcedure{}
\end{algorithmic}
\end{algorithm}

\begin{remark}[Remark on Scaling Tests]
\label{sec:scaling-remark}
Strong scaling tests (where the global sizes $\numparam$ and $\numdata$ are fixed) can be performed by selecting the optimal grid shape using~\cref{alg:grid} for each value of $p$ that is tested. Note however, that it is not possible to keep the local sizes $\procrowblock$ and $\proccolblock$ fixed in addition to $l = \log_{10} \qty(\numdata / \numparam)$ as the number of processors $p$ is increased. The exception to this is the case where there is a square grid ($\numprocrows = \numproccols$), and $p = n^2, n\in \N$. However, a square grid is not the optimal grid shape for any nonzero value of $l$, so performing weak scaling tests on a square grid is not useful for the majority of cases. Instead, if the value of $l$ is allowed to change with $p$ but the \textit{sign} of $l$ remains constant, it is possible to perform weak scaling tests. For this purpose, consider the modified asymptotic communication cost function in~\cref{eq:modcost} written in terms of the local sizes  $\procrowblock$ and $\proccolblock$ (which are now assumed constant)
\begin{align}
    \tilde{C}(\numprocrows) = \proccolblock \log \numprocrows + \procrowblock \log \frac{p}{\numprocrows} = \proccolblock\qty(\log \numprocrows + L\log \frac{p}{\numprocrows}) = \proccolblock\qty(L\log p + (1-L)\log \numprocrows),
\end{align}
where $L = \procrowblock/\proccolblock$. As before, the overall constant scale factor $\proccolblock$ and shift factor $L\log p$ are removed to get the modified asymptotic communication cost for constant local sizes
\begin{align}
    \tilde{C}_{\text{loc}}(\numprocrows) := (1-L)\log \numprocrows.
\end{align}
Now, there are three cases. If $L = 1$ ($\procrowblock = \proccolblock$), the cost is independent of $\numprocrows$. Thus, any grid shape can be used for weak scaling tests. If $L > 1$ ($\procrowblock > \proccolblock$), the cost is a monotonically decreasing function of $\numprocrows$. Thus, a grid shape of $p \times 1$ should be used for weak scaling tests. If $L < 1$, the cost is a monotonically increasing function of $\numprocrows$. Thus, a grid shape of $1\times p$ should be used for weak scaling tests. For most practical applications, $L$ is less than 1.
\end{remark}

%% file: numerical_results.tex
The scalability and performance of the matvec algorithm were tested first on a single GPU, and then on multiple GPUs. For all GPU tests, NVIDIA A100 40~GB GPU was used, and the runtimes of the different steps in the matvec algorithm (see~\cref{tab:ComputationalCost}) were calculated. For single GPU testing, matvecs with $\blocktoep$ and $\blocktoep^*$ were used for different values of $\procrowblock$, $\proccolblock$, and $\numtime$. For multi-GPU testing, matvecs with $\blocktoep$ and $\blocktoep^*$ were tested for different numbers of GPUs $p$. Furthermore, the CUDA kernels used in the algorithm were profiled with the \texttt{NSight Compute} software\footnote{\url{https://developer.nvidia.com/nsight-compute}} to determine their performance. For all tests, matvecs were computed between fixed, predetermined matrices and vectors. 
\revstart 
These matrices and vectors do not correspond to any specific autonomous dynamical system of interest; they are synthetically generated so that the output vector from the matvec has an analytical solution. This allows the correctness of the matvec to be easily verified without explicitly forming the full, dense $\blocktoep$ matrix, which would be intractable for large problem sizes.   
\revend

\revend
Our implementation is available at \url{https://github.com/s769/FFTMatvec}.

For the multi-GPU test, the algorithm was tested on up to 16 GPU nodes of the \emph{Lonestar6} supercomputer at the Texas Advanced Computing Center (TACC). Each GPU node has 3 NVIDIA A100 40~GB GPUs.\footnote{\url{https://docs.tacc.utexas.edu/hpc/lonestar6/\#system-gpu}} These nodes were used to determine the weak and strong scalability of our algorithm. Times for the setup phase of the algorithm are not reported, since these are one-time costs.

\subsection{Single-GPU Performance}\label{sec:SingleGPU}

\begin{figure}[htbp]
    \centering
    \resizebox{\columnwidth}{!}{\input{images/new-single-gpu-scaling.tex}}
    \caption{Total (left) and normalized (right) single-GPU scaling results for the full matvec showing the breakdown of times among different steps of the local matvec. Here, $\numtime = 2{,}000$ and $\numdata = 7$ throughout. The SBGEMV comprises the majority of the total matvec time relative to the other steps. The total computation time per element tends to a constant as the input sizes are increased.}\label{fig:SingleGPUScaling}
\end{figure}
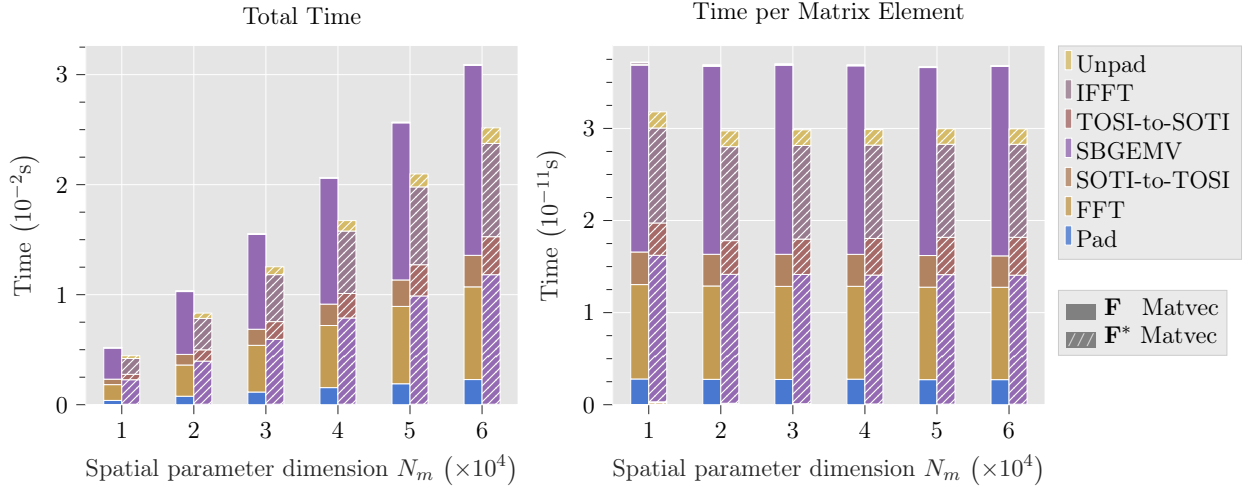

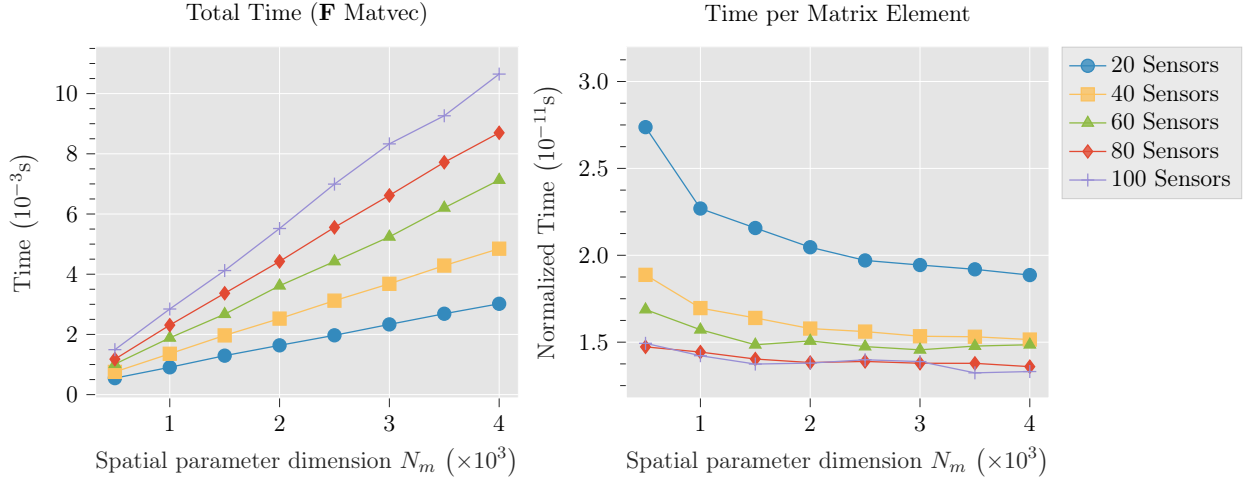
\begin{figure}[htbp]
  \resizebox{\columnwidth}{!}{\input{images/new-single-gpu-scaling-line.tex}}
  \caption{Total (left) and normalized (right) single-GPU scaling results for the full matvec with varying values of $\numparam$ and $\numdata$. Here, $\numtime = 2{,}000$ throughout. Linear scaling with respect to both $\numparam$ and $\numdata$ is observed. The total computation time per element tends to a constant as the input sizes are increased.}\label{fig:SingleGPUScalingLine}
\end{figure}

\Cref{fig:SingleGPUScaling} shows the scaling results for a single GPU---aggregate and normalized by the number of local matrix elements ($2\numtime\procrowblock\proccolblock$). The number of time steps $\numtime$ for this plot is fixed at $2{,}000$, and the spatial data dimension (number of sensors) $\numdata$ is fixed at 7. Each bar represents a different value of the spatial parameter dimension (number of sources) $\numparam$. The rightmost bar represents a total array size close to the memory limit of the single GPU.\@ Each bar is separated into sections showing the times of the different portions of the local matvec. These figures show that the matrix application (SBGEMV) comprises the most time. This is to be expected since it is of size $\numtime\procrowblock\proccolblock$. Also note from the normalized times in \cref{fig:SingleGPUScaling} that as the total number of local matrix elements increases, the total cost per element remains roughly the same. \Cref{fig:SingleGPUScalingLine} shows the total and normalized local matvec times on a single GPU for varying values of $\numparam$ and $\numdata$, with $\numtime$ fixed at $2{,}000$. The expected linear scaling is observed with respect to both $\numparam$ and $\numdata$, and we find that the total computation time per element remains roughly the same as $\numparam$ and $\numdata$ increase.

\begin{figure}[htbp]
  \resizebox{\columnwidth}{!}{\input{images/new-single-gpu-time-scaling-line.tex}}
  \caption{Total (left) and normalized (right) single-GPU scaling results for the full matvec with varying values of $\numtime$. Here, $\numparam = 800$ and $\numdata = 100$ throughout. Approximate linear scaling is observed. The total computation time per element tends slightly decreases as $\numtime$ is increased.
  }\label{fig:SingleGPUTimeScalingLine}
\end{figure}
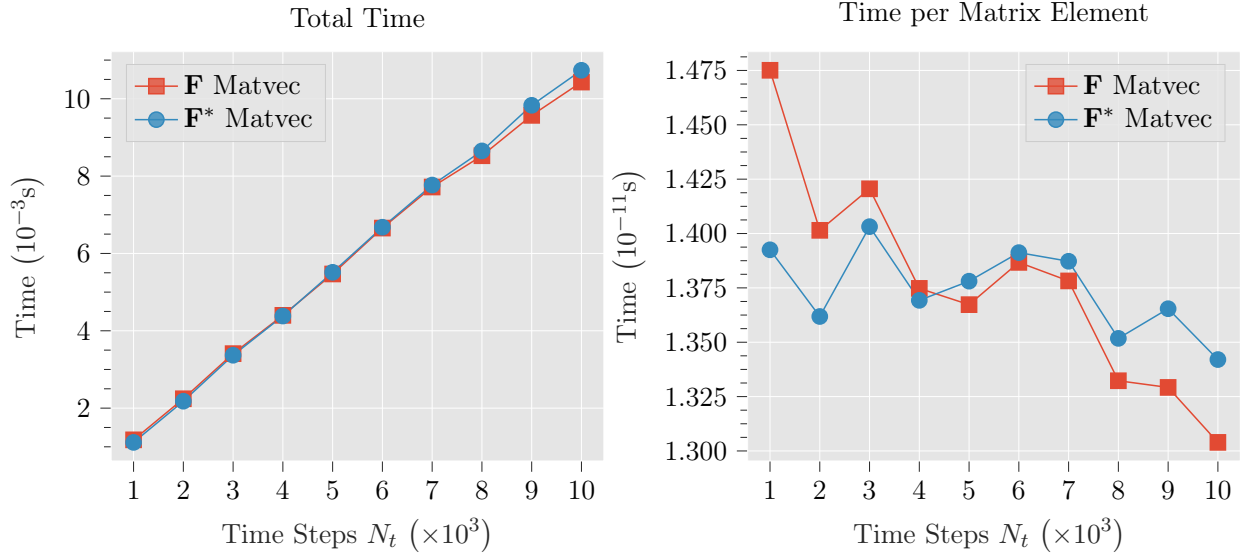

\Cref{fig:SingleGPUTimeScalingLine} shows the total and normalized local matvec times on a single GPU for varying values of $\numtime$, with $\numparam = 800$ and $\numdata = 100$ fixed. According to the computational costs in~\cref{tab:ComputationalCost}, the total matvec time should scale asymptotically as $\numtime \log \numtime$, but practically scales as $\numtime$ for $\numtime < \numparam \numdata\ (=\procrowblock\proccolblock \text{ for a single GPU})$. Approximate linear scaling is observed with respect to $\numtime$, and the total computation time per element slightly decreases as $\numtime$ increases.

\begin{figure}[htpb]
    \centering
    \input{images/new_roofline.tex}
    \caption{Roofline plots for the main kernels used in the matvec computation. All major kernels are operating at close to peak performance. Some kernels are not operating at peak performance due to small input size (they operate on a vector of size $\procrowblock \ll \proccolblock$); these kernels comprise less than 1\% of the total runtime.}\label{fig:Roofline}
\end{figure}
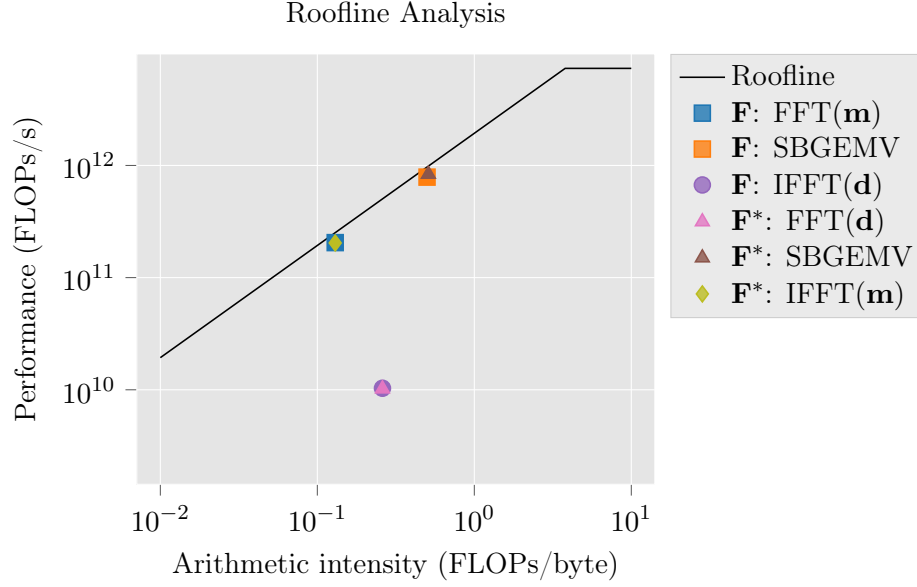

In~\cref{fig:Roofline}, the roofline plots (obtained from \texttt{NSight Compute}) are shown. 
All major kernels (ReIndex not shown on the roofline plot) achieve $75$--$85$\% of the peak memory bandwidth. Note that the apparently poor performance of some kernels (FFT/IFFT$\qty(\datavec)$) is purely due to the small input size passed to them; these kernels account for less than 1\% of the total runtime and do not pose an issue for overall performance. Roofline plots for the padding and index reordering kernels are not included since they do not have any arithmetic operations. However, these kernels are also memory bound and achieve $75$--$85$\% of the peak memory bandwidth. The major computational kernel is the SBGEMV (see~\cref{fig:SingleGPUScaling}). The theoretical FLOP count for this operation (for a single time index) is $8\procrowblock\proccolblock$---the $8$ is for 1 complex add and 1 complex multiply. The number of bytes used for the operation (at a single time index) is $16(\procrowblock\proccolblock + \proccolblock + \procrowblock)$---16 bytes per complex number in the matrix and vectors. Thus, the theoretical arithmetic intensity is $\procrowblock\proccolblock/(2(\procrowblock\proccolblock + \proccolblock + \procrowblock))$. Assuming $\proccolblock \ll \proccolblock \procrowblock$ and $\procrowblock \ll \proccolblock\procrowblock$, the theoretical arithmetic intensity is $0.5$---which exactly matches the observed arithmetic intensity of the SBGEMV kernels. The analysis is the same for both the $\blocktoep$ and $\blocktoep^*$ matvecs.

\subsection{Multi-GPU Performance}\label{sec:MultiGPU}

\begin{figure}[htbp]
  \centering
  \resizebox{\columnwidth}{!}{\input{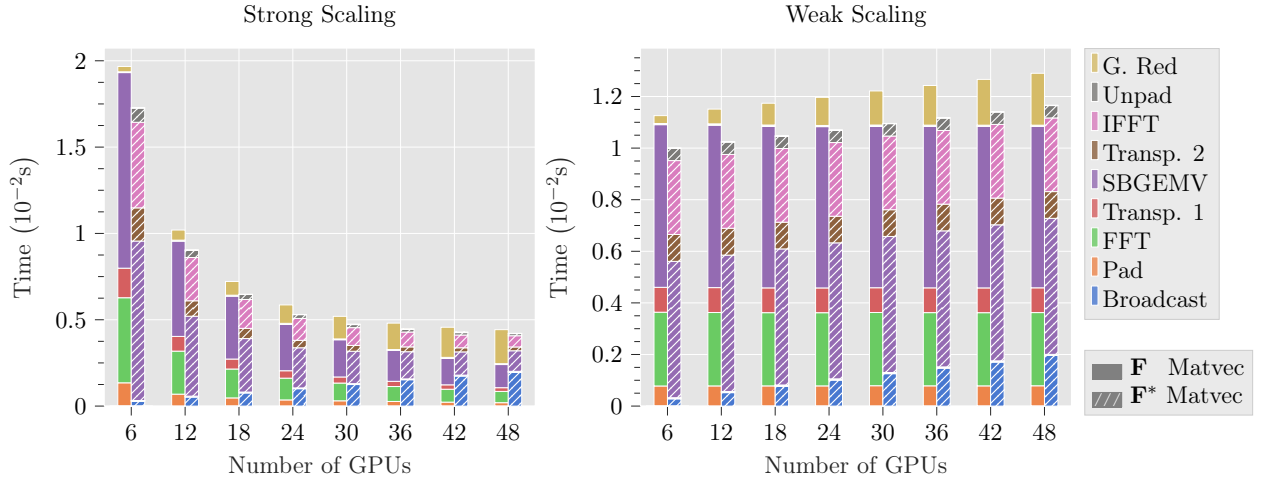}}
    \caption{Strong (left) and weak (right) scaling results for the $\blocktoep$ and $\blocktoep^*$ matvecs on up to 48 NVIDIA A100 40~GB GPUs on TACC's \emph{Lonestar6} supercomputer. The GPUs are arranged in a $1 \times p$ grid, where $p$ is the total number of GPUs.}\label{fig:StrongWeakScalingNoSetup}
\end{figure}

\Cref{fig:StrongWeakScalingNoSetup} shows the strong and weak scaling results split over the different steps of the algorithm. As mentioned in~\cref{sec:scaling-remark}, a $1\times p$ grid of GPUs was used for weak scaling, with local sizes $\procrowblock = 10$, $\numtime = 2{,}000$, $\proccolblock = 20{,}000$. A $1\times p$ grid was also used for strong scaling, and the global problem size was set to the maximum size that could fit on 6~GPUs ($\procrowblock = 10$, $\numtime = 2{,}000$, $\proccolblock = 210{,}000$). \Cref{fig:StrongWeakScalingNoSetup} shows that the computation costs remain constant as the number of GPUs is increased, and the communication costs increase at a sub-logarithmic rate. This is likely due to differences in on-node vs.\ off-node communication costs as explained in~\cref{sec:gpu-grid}. Communication quickly begins to dominate the total runtime for strong scaling experiments.

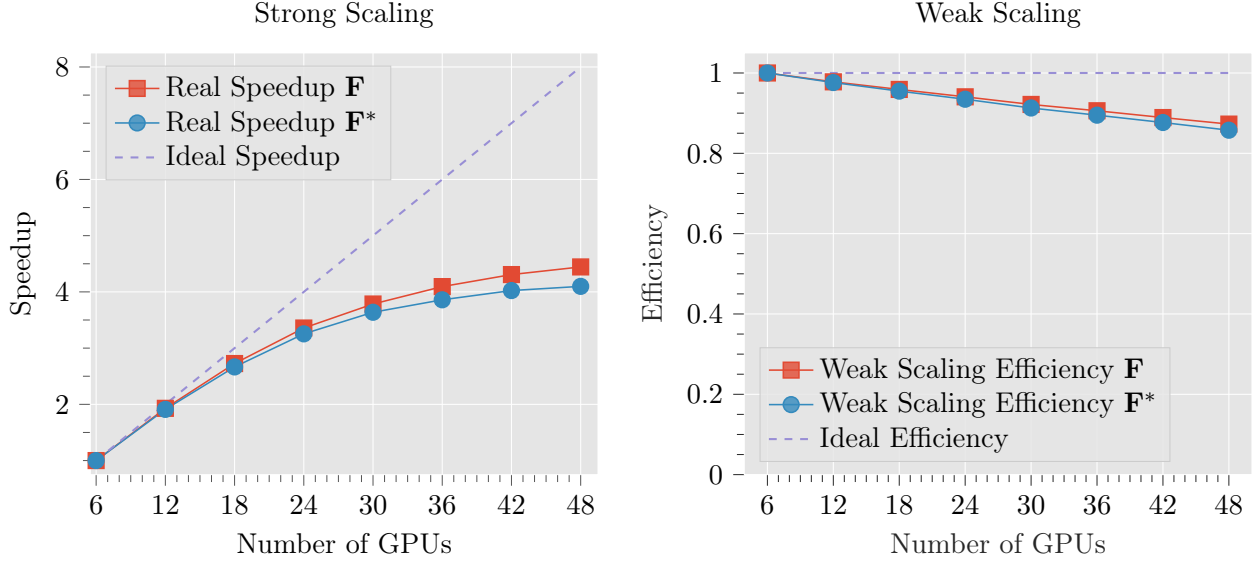
\begin{figure}[htbp]
  \centering
  \resizebox{\columnwidth}{!}{\input{images/new_weak_scaling_full_ls6}}
  \caption{Strong (left) and weak (right) scaling results for the $\blocktoep$ and $\blocktoep^*$ matvecs on up to 48 NVIDIA A100 40~GB GPUs on TACC's \emph{Lonestar6} supercomputer. All tests are done on a $1 \times p$ grid, where $p$ is the total number of GPUs. The weak scaling efficiency is controlled by the communication costs, which are growing sub-logarithmically as seen in~\cref{fig:StrongWeakScalingNoSetup}. The drop in strong scaling efficiency is due to dominating communication costs and the fact that smaller problem sizes decrease GPU kernel efficiency. This effect is also observed in~\cref{fig:SingleGPUScalingLine}.}\label{fig:StrongWeakScalingFull}
\end{figure}

\Cref{fig:StrongWeakScalingFull} shows the results of the strong and weak scaling tests. All tests are performed on a $1\times p$ grid, where $p$ is the total number of GPUs. This is the optimal grid shape according to~\cref{alg:grid} for all cases. The weak scaling efficiency is controlled by the communication costs, which are growing sub-logarithmically as seen in~\cref{fig:StrongWeakScalingNoSetup}. The loss of strong scaling speedup can be explained by increased communication costs and the fact that smaller problem sizes decrease GPU kernel efficiency. This effect is also observed in~\cref{fig:SingleGPUScalingLine}, where decreasing the problem size increases the computational time per element of the matrix.

It is important to note that the goal of this algorithm is to enable real-time solutions of inverse problems. From that perspective, what matters is the overall time to solution rather than scalability. In that regard, the matvec portion of this algorithm provides a significant speedup over conventional methods as mentioned in~\cref{sec:InverseProblemsInvolvingShiftInvariantSystems}. A detailed comparison of the matvec runtimes for this vs.~conventional algorithms will be presented in a follow-up paper. As a preview, note that for one such example, this algorithm provides a $\sim 10{,}000\times$ speedup over conventional methods.\footnote{The model used for this comparison is the acoustic-gravity wave equations that model tsunami dynamics~\cite{lotto2015tsunami}.}

\begin{figure}[htbp]
  \centering
  \input{images/grid_test}
  \caption{Communication times vs. number of rows in a processor grid with 80 GPUs for different values of $l = \log_{10}(\numdata/\numparam)$. Minimum times for each $l$ value are scaled to 1. The observed minimal number of rows for each $l$ value matches the theoretical minimal number of rows from~\cref{alg:grid}.}\label{fig:GridTest}
\end{figure}
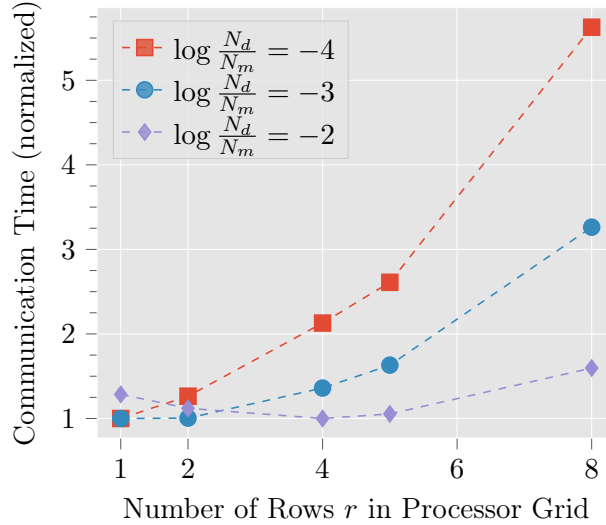

To compare the observed optimal grid shape for a given problem configuration to the theoretical optimal grid shape from~\cref{alg:grid}, several grid shapes for 80 GPUs were tested with different values of $l = \log_{10}(\numdata/\numparam)$. These tests were carried out on TACC's \emph{Frontera} GPU nodes (with four NVIDIA Quadro RTX 5000 GPUs per node).\footnote{\url{https://docs.tacc.utexas.edu/hpc/frontera/\#system-gpu}} The results are shown in~\cref{fig:GridTest}, with minimum times for each value of $l$ scaled to 1. The observed minimal number of rows matches the theoretical number for each value of $l$. Only values of $l \leq -2$ are considered since this is most often the case in practice (many more parameters than data).

%% file: images/new-single-gpu-scaling.tex
\begin{tikzpicture}

  \definecolor{crimson2143940}{RGB}{214,39,40}
  \definecolor{darkorange25512714}{RGB}{255,127,14}
  \definecolor{darkslategray38}{RGB}{38,38,38}
  \definecolor{darkturquoise23190207}{RGB}{23,190,207}
  \definecolor{forestgreen4416044}{RGB}{44,160,44}
  \definecolor{goldenrod18818934}{RGB}{188,189,34}
  \definecolor{gray127}{RGB}{127,127,127}
  \definecolor{lavender234234242}{RGB}{234,234,242}
  \definecolor{lightblue158218229}{RGB}{158,218,229}
  \definecolor{lightgray204}{RGB}{204,204,204}
  \definecolor{mediumpurple148103189}{RGB}{148,103,189}
  \definecolor{orchid227119194}{RGB}{227,119,194}
  \definecolor{sienna1408675}{RGB}{140,86,75}
  \definecolor{steelblue31119180}{RGB}{31,119,180}
  \definecolor{gainsboro229}{RGB}{229,229,229}

\definecolor{darkkhaki213187103}{RGB}{213,187,103}
\definecolor{gray}{RGB}{128,128,128}
\definecolor{gray152122143}{RGB}{152,122,143}
\definecolor{indianred166106104}{RGB}{166,106,104}
\definecolor{mediumpurple148107178}{RGB}{148,107,178}
\definecolor{peru17813196}{RGB}{178,131,96}
\definecolor{peru19615582}{RGB}{196,155,82}
\definecolor{royalblue72120208}{RGB}{72,120,208}

\begin{axis}[
title={Total Time},
axis background/.style={fill=gainsboro229},
name=ax1,
axis line style={white},
legend cell align={left},
legend style={
  fill opacity=0.8,
  draw opacity=1,
  text opacity=1,
  at={(0.03,0.97)},
  anchor=north west,
  draw=lightgray204,
  fill=gainsboro229
},
legend pos=outer north east,
reverse legend,
minor y tick num=3,
y tick label style={/pgf/number format/fixed},
tick align=outside,
tick pos=left,
x grid style={white},
xlabel=\textcolor{darkslategray38}{Spatial parameter dimension $\numparam$ \(\displaystyle \qty(\times 10^4)\)},
xmajorgrids,
xmin=-0.5, xmax=5.5,
xtick style={color=darkslategray38},
xtick={0,1,2,3,4,5},
xticklabels={1,2,3,4,5,6},
y grid style={white},
ylabel=\textcolor{darkslategray38}{Time $\qty(10^{-2}\text{s})$},
scaled ticks=false,
ymajorgrids,
ymin=0, ymax=0.0326467365,
ytick={0,0.01,0.02,0.03,0.04},
yticklabels={0,1,2,3,4},
ytick style={color=darkslategray38}
]
\draw[draw=white,fill=royalblue72120208,line width=0pt] (axis cs:-0.25,0) rectangle (axis cs:0,0.00039115);
\addlegendimage{ybar,ybar legend,draw=white,fill=royalblue72120208,line width=0pt}
\addlegendentry{Pad}

\draw[draw=white,fill=royalblue72120208,line width=0pt] (axis cs:0.75,0) rectangle (axis cs:1,0.00077144);
\draw[draw=white,fill=royalblue72120208,line width=0pt] (axis cs:1.75,0) rectangle (axis cs:2,0.0011525);
\draw[draw=white,fill=royalblue72120208,line width=0pt] (axis cs:2.75,0) rectangle (axis cs:3,0.00155296);
\draw[draw=white,fill=royalblue72120208,line width=0pt] (axis cs:3.75,0) rectangle (axis cs:4,0.00190466);
\draw[draw=white,fill=royalblue72120208,line width=0pt] (axis cs:4.75,0) rectangle (axis cs:5,0.0022861);
\draw[draw=white,fill=peru19615582,line width=0pt] (axis cs:-0.25,0.00039115) rectangle (axis cs:0,0.00182616);
\addlegendimage{ybar,ybar legend,draw=white,fill=peru19615582,line width=0pt}
\addlegendentry{FFT}

\draw[draw=white,fill=peru19615582,line width=0pt] (axis cs:0.75,0.00077144) rectangle (axis cs:1,0.00361);
\draw[draw=white,fill=peru19615582,line width=0pt] (axis cs:1.75,0.0011525) rectangle (axis cs:2,0.0053964);
\draw[draw=white,fill=peru19615582,line width=0pt] (axis cs:2.75,0.00155296) rectangle (axis cs:3,0.00719859);
\draw[draw=white,fill=peru19615582,line width=0pt] (axis cs:3.75,0.00190466) rectangle (axis cs:4,0.00893424);
\draw[draw=white,fill=peru19615582,line width=0pt] (axis cs:4.75,0.0022861) rectangle (axis cs:5,0.0107121);
\draw[draw=white,fill=peru17813196,line width=0pt] (axis cs:-0.25,0.00182616) rectangle (axis cs:0,0.00232158);
\addlegendimage{ybar,ybar legend,draw=white,fill=peru17813196,line width=0pt}
\addlegendentry{SOTI-to-TOSI}

\draw[draw=white,fill=peru17813196,line width=0pt] (axis cs:0.75,0.00361) rectangle (axis cs:1,0.00457333);
\draw[draw=white,fill=peru17813196,line width=0pt] (axis cs:1.75,0.0053964) rectangle (axis cs:2,0.0068591);
\draw[draw=white,fill=peru17813196,line width=0pt] (axis cs:2.75,0.00719859) rectangle (axis cs:3,0.00914307);
\draw[draw=white,fill=peru17813196,line width=0pt] (axis cs:3.75,0.00893424) rectangle (axis cs:4,0.01133919);
\draw[draw=white,fill=peru17813196,line width=0pt] (axis cs:4.75,0.0107121) rectangle (axis cs:5,0.0135757);
\draw[draw=white,fill=mediumpurple148107178,line width=0pt] (axis cs:-0.25,0.00232158) rectangle (axis cs:0,0.00516141);
\addlegendimage{ybar,ybar legend,draw=white,fill=mediumpurple148107178,line width=0pt}
\addlegendentry{SBGEMV}

\draw[draw=white,fill=mediumpurple148107178,line width=0pt] (axis cs:0.75,0.00457333) rectangle (axis cs:1,0.01028324);
\draw[draw=white,fill=mediumpurple148107178,line width=0pt] (axis cs:1.75,0.0068591) rectangle (axis cs:2,0.01549095);
\draw[draw=white,fill=mediumpurple148107178,line width=0pt] (axis cs:2.75,0.00914307) rectangle (axis cs:3,0.02060907);
\draw[draw=white,fill=mediumpurple148107178,line width=0pt] (axis cs:3.75,0.01133919) rectangle (axis cs:4,0.02562636);
\draw[draw=white,fill=mediumpurple148107178,line width=0pt] (axis cs:4.75,0.0135757) rectangle (axis cs:5,0.03087194);
\draw[draw=white,fill=indianred166106104,line width=0pt] (axis cs:-0.25,0.00516141) rectangle (axis cs:0,0.0051739);
\addlegendimage{ybar,ybar legend,draw=white,fill=indianred166106104,line width=0pt}
\addlegendentry{TOSI-to-SOTI}

\draw[draw=white,fill=indianred166106104,line width=0pt] (axis cs:0.75,0.01028324) rectangle (axis cs:1,0.01029572);
\draw[draw=white,fill=indianred166106104,line width=0pt] (axis cs:1.75,0.01549095) rectangle (axis cs:2,0.01550274);
\draw[draw=white,fill=indianred166106104,line width=0pt] (axis cs:2.75,0.02060907) rectangle (axis cs:3,0.02062062);
\draw[draw=white,fill=indianred166106104,line width=0pt] (axis cs:3.75,0.02562636) rectangle (axis cs:4,0.02563851);
\draw[draw=white,fill=indianred166106104,line width=0pt] (axis cs:4.75,0.03087194) rectangle (axis cs:5,0.03088482);
\draw[draw=white,fill=gray152122143,line width=0pt] (axis cs:-0.25,0.0051739) rectangle (axis cs:0,0.00519916);
\addlegendimage{ybar,ybar legend,draw=white,fill=gray152122143,line width=0pt}
\addlegendentry{IFFT}

\draw[draw=white,fill=gray152122143,line width=0pt] (axis cs:0.75,0.01029572) rectangle (axis cs:1,0.01032105);
\draw[draw=white,fill=gray152122143,line width=0pt] (axis cs:1.75,0.01550274) rectangle (axis cs:2,0.01552716);
\draw[draw=white,fill=gray152122143,line width=0pt] (axis cs:2.75,0.02062062) rectangle (axis cs:3,0.02064675);
\draw[draw=white,fill=gray152122143,line width=0pt] (axis cs:3.75,0.02563851) rectangle (axis cs:4,0.02566346);
\draw[draw=white,fill=gray152122143,line width=0pt] (axis cs:4.75,0.03088482) rectangle (axis cs:5,0.03091051);
\draw[draw=white,fill=darkkhaki213187103,line width=0pt] (axis cs:-0.25,0.00519916) rectangle (axis cs:0,0.00520961);
\addlegendimage{ybar,ybar legend,draw=white,fill=darkkhaki213187103,line width=0pt}
\addlegendentry{Unpad}

\draw[draw=white,fill=darkkhaki213187103,line width=0pt] (axis cs:0.75,0.01032105) rectangle (axis cs:1,0.01033192);
\draw[draw=white,fill=darkkhaki213187103,line width=0pt] (axis cs:1.75,0.01552716) rectangle (axis cs:2,0.01553776);
\draw[draw=white,fill=darkkhaki213187103,line width=0pt] (axis cs:2.75,0.02064675) rectangle (axis cs:3,0.02065766);
\draw[draw=white,fill=darkkhaki213187103,line width=0pt] (axis cs:3.75,0.02566346) rectangle (axis cs:4,0.02567432);
\draw[draw=white,fill=darkkhaki213187103,line width=0pt] (axis cs:4.75,0.03091051) rectangle (axis cs:5,0.0309213);
\draw[draw=white,fill=royalblue72120208,line width=0pt,postaction={pattern=north east lines, pattern color=white}] (axis cs:0,0) rectangle (axis cs:0.25,1.136e-05);
\draw[draw=white,fill=royalblue72120208,line width=0pt,postaction={pattern=north east lines, pattern color=white}] (axis cs:1,0) rectangle (axis cs:1.25,1.155e-05);
\draw[draw=white,fill=royalblue72120208,line width=0pt,postaction={pattern=north east lines, pattern color=white}] (axis cs:2,0) rectangle (axis cs:2.25,1.142e-05);
\draw[draw=white,fill=royalblue72120208,line width=0pt,postaction={pattern=north east lines, pattern color=white}] (axis cs:3,0) rectangle (axis cs:3.25,1.154e-05);
\draw[draw=white,fill=royalblue72120208,line width=0pt,postaction={pattern=north east lines, pattern color=white}] (axis cs:4,0) rectangle (axis cs:4.25,1.173e-05);
\draw[draw=white,fill=royalblue72120208,line width=0pt,postaction={pattern=north east lines, pattern color=white}] (axis cs:5,0) rectangle (axis cs:5.25,1.172e-05);
\draw[draw=white,fill=peru19615582,line width=0pt,postaction={pattern=north east lines, pattern color=white}] (axis cs:0,1.136e-05) rectangle (axis cs:0.25,3.36e-05);
\draw[draw=white,fill=peru19615582,line width=0pt,postaction={pattern=north east lines, pattern color=white}] (axis cs:1,1.155e-05) rectangle (axis cs:1.25,3.4e-05);
\draw[draw=white,fill=peru19615582,line width=0pt,postaction={pattern=north east lines, pattern color=white}] (axis cs:2,1.142e-05) rectangle (axis cs:2.25,3.345e-05);
\draw[draw=white,fill=peru19615582,line width=0pt,postaction={pattern=north east lines, pattern color=white}] (axis cs:3,1.154e-05) rectangle (axis cs:3.25,3.387e-05);
\draw[draw=white,fill=peru19615582,line width=0pt,postaction={pattern=north east lines, pattern color=white}] (axis cs:4,1.173e-05) rectangle (axis cs:4.25,3.418e-05);
\draw[draw=white,fill=peru19615582,line width=0pt,postaction={pattern=north east lines, pattern color=white}] (axis cs:5,1.172e-05) rectangle (axis cs:5.25,3.479e-05);
\draw[draw=white,fill=peru17813196,line width=0pt,postaction={pattern=north east lines, pattern color=white}] (axis cs:0,3.36e-05) rectangle (axis cs:0.25,4.414e-05);
\draw[draw=white,fill=peru17813196,line width=0pt,postaction={pattern=north east lines, pattern color=white}] (axis cs:1,3.4e-05) rectangle (axis cs:1.25,4.478e-05);
\draw[draw=white,fill=peru17813196,line width=0pt,postaction={pattern=north east lines, pattern color=white}] (axis cs:2,3.345e-05) rectangle (axis cs:2.25,4.414e-05);
\draw[draw=white,fill=peru17813196,line width=0pt,postaction={pattern=north east lines, pattern color=white}] (axis cs:3,3.387e-05) rectangle (axis cs:3.25,4.454e-05);
\draw[draw=white,fill=peru17813196,line width=0pt,postaction={pattern=north east lines, pattern color=white}] (axis cs:4,3.418e-05) rectangle (axis cs:4.25,4.493e-05);
\draw[draw=white,fill=peru17813196,line width=0pt,postaction={pattern=north east lines, pattern color=white}] (axis cs:5,3.479e-05) rectangle (axis cs:5.25,4.572e-05);
\draw[draw=white,fill=mediumpurple148107178,line width=0pt,postaction={pattern=north east lines, pattern color=white}] (axis cs:0,4.414e-05) rectangle (axis cs:0.25,0.00227144);
\draw[draw=white,fill=mediumpurple148107178,line width=0pt,postaction={pattern=north east lines, pattern color=white}] (axis cs:1,4.478e-05) rectangle (axis cs:1.25,0.00396369);
\draw[draw=white,fill=mediumpurple148107178,line width=0pt,postaction={pattern=north east lines, pattern color=white}] (axis cs:2,4.414e-05) rectangle (axis cs:2.25,0.00593985);
\draw[draw=white,fill=mediumpurple148107178,line width=0pt,postaction={pattern=north east lines, pattern color=white}] (axis cs:3,4.454e-05) rectangle (axis cs:3.25,0.00788347);
\draw[draw=white,fill=mediumpurple148107178,line width=0pt,postaction={pattern=north east lines, pattern color=white}] (axis cs:4,4.493e-05) rectangle (axis cs:4.25,0.00988607);
\draw[draw=white,fill=mediumpurple148107178,line width=0pt,postaction={pattern=north east lines, pattern color=white}] (axis cs:5,4.572e-05) rectangle (axis cs:5.25,0.01181832);
\draw[draw=white,fill=indianred166106104,line width=0pt,postaction={pattern=north east lines, pattern color=white}] (axis cs:0,0.00227144) rectangle (axis cs:0.25,0.00276264);
\draw[draw=white,fill=indianred166106104,line width=0pt,postaction={pattern=north east lines, pattern color=white}] (axis cs:1,0.00396369) rectangle (axis cs:1.25,0.00498855);
\draw[draw=white,fill=indianred166106104,line width=0pt,postaction={pattern=north east lines, pattern color=white}] (axis cs:2,0.00593985) rectangle (axis cs:2.25,0.00755245);
\draw[draw=white,fill=indianred166106104,line width=0pt,postaction={pattern=north east lines, pattern color=white}] (axis cs:3,0.00788347) rectangle (axis cs:3.25,0.01010978);
\draw[draw=white,fill=indianred166106104,line width=0pt,postaction={pattern=north east lines, pattern color=white}] (axis cs:4,0.00988607) rectangle (axis cs:4.25,0.01272417);
\draw[draw=white,fill=indianred166106104,line width=0pt,postaction={pattern=north east lines, pattern color=white}] (axis cs:5,0.01181832) rectangle (axis cs:5.25,0.01527089);
\draw[draw=white,fill=gray152122143,line width=0pt,postaction={pattern=north east lines, pattern color=white}] (axis cs:0,0.00276264) rectangle (axis cs:0.25,0.00420578);
\draw[draw=white,fill=gray152122143,line width=0pt,postaction={pattern=north east lines, pattern color=white}] (axis cs:1,0.00498855) rectangle (axis cs:1.25,0.00784572);
\draw[draw=white,fill=gray152122143,line width=0pt,postaction={pattern=north east lines, pattern color=white}] (axis cs:2,0.00755245) rectangle (axis cs:2.25,0.01182226);
\draw[draw=white,fill=gray152122143,line width=0pt,postaction={pattern=north east lines, pattern color=white}] (axis cs:3,0.01010978) rectangle (axis cs:3.25,0.01578526);
\draw[draw=white,fill=gray152122143,line width=0pt,postaction={pattern=north east lines, pattern color=white}] (axis cs:4,0.01272417) rectangle (axis cs:4.25,0.0197961);
\draw[draw=white,fill=gray152122143,line width=0pt,postaction={pattern=north east lines, pattern color=white}] (axis cs:5,0.01527089) rectangle (axis cs:5.25,0.02374931);
\draw[draw=white,fill=darkkhaki213187103,line width=0pt,postaction={pattern=north east lines, pattern color=white}] (axis cs:0,0.00420578) rectangle (axis cs:0.25,0.00445178);
\draw[draw=white,fill=darkkhaki213187103,line width=0pt,postaction={pattern=north east lines, pattern color=white}] (axis cs:1,0.00784572) rectangle (axis cs:1.25,0.00832602);
\draw[draw=white,fill=darkkhaki213187103,line width=0pt,postaction={pattern=north east lines, pattern color=white}] (axis cs:2,0.01182226) rectangle (axis cs:2.25,0.01253639);
\draw[draw=white,fill=darkkhaki213187103,line width=0pt,postaction={pattern=north east lines, pattern color=white}] (axis cs:3,0.01578526) rectangle (axis cs:3.25,0.01673245);
\draw[draw=white,fill=darkkhaki213187103,line width=0pt,postaction={pattern=north east lines, pattern color=white}] (axis cs:4,0.0197961) rectangle (axis cs:4.25,0.02097454);
\draw[draw=white,fill=darkkhaki213187103,line width=0pt,postaction={pattern=north east lines, pattern color=white}] (axis cs:5,0.02374931) rectangle (axis cs:5.25,0.02516124);
\draw[draw=white,fill=gray,very thin] (axis cs:-0.4,0) rectangle (axis cs:0.4,0);
\draw[draw=white,fill=gray,very thin,postaction={pattern=north east lines, pattern color=white}] (axis cs:-0.4,0) rectangle (axis cs:0.4,0);
\legend{};
\end{axis}

\begin{axis}[
  title={Time per Matrix Element},
  at={(ax1.south east)},
  name=ax2,
  xshift=1.5cm,
  hide scale,
  axis background/.style={fill=gainsboro229},
  axis line style={white},
  legend cell align={left},
  legend style={
    fill opacity=0.8,
    draw opacity=1,
    text opacity=1,
    at={(0.09,0.5)},
    anchor=west,
    draw=lightgray204,
    fill=gainsboro229,
    name=leg
  },
  legend pos=outer north east,
  reverse legend,
  tick align=outside,
  minor y tick num=3,
  tick pos=left,
  x grid style={white},
  xlabel=\textcolor{darkslategray38}{Spatial parameter dimension $\numparam$ \(\displaystyle \qty(\times 10^4)\)},
  xmajorgrids,
  xmin=-0.5, xmax=5.5,
xtick style={color=darkslategray38},
xtick={0,1,2,3,4,5},
xticklabels={1,2,3,4,5,6},
  y grid style={white},
  ylabel=\textcolor{darkslategray38}{Time $\qty(10^{-11}\text{s})$},
  ymajorgrids,
  scaled ticks=false,
  ymin=0, ymax=3.9e-11,
  ytick style={color=darkslategray38},
  ytick={0,1e-11,2e-11,3e-11,4e-11},
  yticklabels={0,1,2,3,4}
  ]
  \draw[draw=white,fill=royalblue72120208,line width=0pt] (axis cs:-0.25,0) rectangle (axis cs:0,2.79392857142857e-12);
\addlegendimage{ybar,ybar legend,draw=white,fill=royalblue72120208,line width=0pt}
\addlegendentry{Pad}

\draw[draw=white,fill=royalblue72120208,line width=0pt] (axis cs:0.75,0) rectangle (axis cs:1,2.75514285714286e-12);
\draw[draw=white,fill=royalblue72120208,line width=0pt] (axis cs:1.75,0) rectangle (axis cs:2,2.74404761904762e-12);
\draw[draw=white,fill=royalblue72120208,line width=0pt] (axis cs:2.75,0) rectangle (axis cs:3,2.77314285714286e-12);
\draw[draw=white,fill=royalblue72120208,line width=0pt] (axis cs:3.75,0) rectangle (axis cs:4,2.72094285714286e-12);
\draw[draw=white,fill=royalblue72120208,line width=0pt] (axis cs:4.75,0) rectangle (axis cs:5,2.72154761904762e-12);
\draw[draw=white,fill=peru19615582,line width=0pt] (axis cs:-0.25,2.79392857142857e-12) rectangle (axis cs:0,1.3044e-11);
\addlegendimage{ybar,ybar legend,draw=white,fill=peru19615582,line width=0pt}
\addlegendentry{FFT}

\draw[draw=white,fill=peru19615582,line width=0pt] (axis cs:0.75,2.75514285714286e-12) rectangle (axis cs:1,1.28928571428571e-11);
\draw[draw=white,fill=peru19615582,line width=0pt] (axis cs:1.75,2.74404761904762e-12) rectangle (axis cs:2,1.28485714285714e-11);
\draw[draw=white,fill=peru19615582,line width=0pt] (axis cs:2.75,2.77314285714286e-12) rectangle (axis cs:3,1.2854625e-11);
\draw[draw=white,fill=peru19615582,line width=0pt] (axis cs:3.75,2.72094285714286e-12) rectangle (axis cs:4,1.27632e-11);
\draw[draw=white,fill=peru19615582,line width=0pt] (axis cs:4.75,2.72154761904762e-12) rectangle (axis cs:5,1.27525e-11);
\draw[draw=white,fill=peru17813196,line width=0pt] (axis cs:-0.25,1.3044e-11) rectangle (axis cs:0,1.65827142857143e-11);
\addlegendimage{ybar,ybar legend,draw=white,fill=peru17813196,line width=0pt}
\addlegendentry{SOTI-to-TOSI}

\draw[draw=white,fill=peru17813196,line width=0pt] (axis cs:0.75,1.28928571428571e-11) rectangle (axis cs:1,1.63333214285714e-11);
\draw[draw=white,fill=peru17813196,line width=0pt] (axis cs:1.75,1.28485714285714e-11) rectangle (axis cs:2,1.63311904761905e-11);
\draw[draw=white,fill=peru17813196,line width=0pt] (axis cs:2.75,1.2854625e-11) rectangle (axis cs:3,1.63269107142857e-11);
\draw[draw=white,fill=peru17813196,line width=0pt] (axis cs:3.75,1.27632e-11) rectangle (axis cs:4,1.61988428571429e-11);
\draw[draw=white,fill=peru17813196,line width=0pt] (axis cs:4.75,1.27525e-11) rectangle (axis cs:5,1.61615476190476e-11);
\draw[draw=white,fill=mediumpurple148107178,line width=0pt] (axis cs:-0.25,1.65827142857143e-11) rectangle (axis cs:0,3.68672142857143e-11);
\addlegendimage{ybar,ybar legend,draw=white,fill=mediumpurple148107178,line width=0pt}
\addlegendentry{SBGEMV}

\draw[draw=white,fill=mediumpurple148107178,line width=0pt] (axis cs:0.75,1.63333214285714e-11) rectangle (axis cs:1,3.67258571428571e-11);
\draw[draw=white,fill=mediumpurple148107178,line width=0pt] (axis cs:1.75,1.63311904761905e-11) rectangle (axis cs:2,3.68832142857143e-11);
\draw[draw=white,fill=mediumpurple148107178,line width=0pt] (axis cs:2.75,1.63269107142857e-11) rectangle (axis cs:3,3.68019107142857e-11);
\draw[draw=white,fill=mediumpurple148107178,line width=0pt] (axis cs:3.75,1.61988428571429e-11) rectangle (axis cs:4,3.66090857142857e-11);
\draw[draw=white,fill=mediumpurple148107178,line width=0pt] (axis cs:4.75,1.61615476190476e-11) rectangle (axis cs:5,3.67523095238095e-11);
\draw[draw=white,fill=indianred166106104,line width=0pt] (axis cs:-0.25,3.68672142857143e-11) rectangle (axis cs:0,3.69564285714286e-11);
\addlegendimage{ybar,ybar legend,draw=white,fill=indianred166106104,line width=0pt}
\addlegendentry{TOSI-to-SOTI}

\draw[draw=white,fill=indianred166106104,line width=0pt] (axis cs:0.75,3.67258571428571e-11) rectangle (axis cs:1,3.67704285714286e-11);
\draw[draw=white,fill=indianred166106104,line width=0pt] (axis cs:1.75,3.68832142857143e-11) rectangle (axis cs:2,3.69112857142857e-11);
\draw[draw=white,fill=indianred166106104,line width=0pt] (axis cs:2.75,3.68019107142857e-11) rectangle (axis cs:3,3.68225357142857e-11);
\draw[draw=white,fill=indianred166106104,line width=0pt] (axis cs:3.75,3.66090857142857e-11) rectangle (axis cs:4,3.66264428571429e-11);
\draw[draw=white,fill=indianred166106104,line width=0pt] (axis cs:4.75,3.67523095238095e-11) rectangle (axis cs:5,3.67676428571429e-11);
\draw[draw=white,fill=gray152122143,line width=0pt] (axis cs:-0.25,3.69564285714286e-11) rectangle (axis cs:0,3.71368571428571e-11);
\addlegendimage{ybar,ybar legend,draw=white,fill=gray152122143,line width=0pt}
\addlegendentry{IFFT}

\draw[draw=white,fill=gray152122143,line width=0pt] (axis cs:0.75,3.67704285714286e-11) rectangle (axis cs:1,3.68608928571429e-11);
\draw[draw=white,fill=gray152122143,line width=0pt] (axis cs:1.75,3.69112857142857e-11) rectangle (axis cs:2,3.69694285714286e-11);
\draw[draw=white,fill=gray152122143,line width=0pt] (axis cs:2.75,3.68225357142857e-11) rectangle (axis cs:3,3.68691964285714e-11);
\draw[draw=white,fill=gray152122143,line width=0pt] (axis cs:3.75,3.66264428571429e-11) rectangle (axis cs:4,3.66620857142857e-11);
\draw[draw=white,fill=gray152122143,line width=0pt] (axis cs:4.75,3.67676428571429e-11) rectangle (axis cs:5,3.67982261904762e-11);
\draw[draw=white,fill=darkkhaki213187103,line width=0pt] (axis cs:-0.25,3.71368571428571e-11) rectangle (axis cs:0,3.72115e-11);
\addlegendimage{ybar,ybar legend,draw=white,fill=darkkhaki213187103,line width=0pt}
\addlegendentry{Unpad}

\draw[draw=white,fill=darkkhaki213187103,line width=0pt] (axis cs:0.75,3.68608928571429e-11) rectangle (axis cs:1,3.68997142857143e-11);
\draw[draw=white,fill=darkkhaki213187103,line width=0pt] (axis cs:1.75,3.69694285714286e-11) rectangle (axis cs:2,3.69946666666667e-11);
\draw[draw=white,fill=darkkhaki213187103,line width=0pt] (axis cs:2.75,3.68691964285714e-11) rectangle (axis cs:3,3.68886785714286e-11);
\draw[draw=white,fill=darkkhaki213187103,line width=0pt] (axis cs:3.75,3.66620857142857e-11) rectangle (axis cs:4,3.66776e-11);
\draw[draw=white,fill=darkkhaki213187103,line width=0pt] (axis cs:4.75,3.67982261904762e-11) rectangle (axis cs:5,3.68110714285714e-11);
\draw[draw=white,fill=royalblue72120208,line width=0pt,postaction={pattern=north east lines, pattern color=white}] (axis cs:0,0) rectangle (axis cs:0.25,8.11428571428571e-14);
\draw[draw=white,fill=royalblue72120208,line width=0pt,postaction={pattern=north east lines, pattern color=white}] (axis cs:1,0) rectangle (axis cs:1.25,4.125e-14);
\draw[draw=white,fill=royalblue72120208,line width=0pt,postaction={pattern=north east lines, pattern color=white}] (axis cs:2,0) rectangle (axis cs:2.25,2.71904761904762e-14);
\draw[draw=white,fill=royalblue72120208,line width=0pt,postaction={pattern=north east lines, pattern color=white}] (axis cs:3,0) rectangle (axis cs:3.25,2.06071428571429e-14);
\draw[draw=white,fill=royalblue72120208,line width=0pt,postaction={pattern=north east lines, pattern color=white}] (axis cs:4,0) rectangle (axis cs:4.25,1.67571428571429e-14);
\draw[draw=white,fill=royalblue72120208,line width=0pt,postaction={pattern=north east lines, pattern color=white}] (axis cs:5,0) rectangle (axis cs:5.25,1.3952380952381e-14);
\draw[draw=white,fill=peru19615582,line width=0pt,postaction={pattern=north east lines, pattern color=white}] (axis cs:0,8.11428571428571e-14) rectangle (axis cs:0.25,2.4e-13);
\draw[draw=white,fill=peru19615582,line width=0pt,postaction={pattern=north east lines, pattern color=white}] (axis cs:1,4.125e-14) rectangle (axis cs:1.25,1.21428571428571e-13);
\draw[draw=white,fill=peru19615582,line width=0pt,postaction={pattern=north east lines, pattern color=white}] (axis cs:2,2.71904761904762e-14) rectangle (axis cs:2.25,7.96428571428571e-14);
\draw[draw=white,fill=peru19615582,line width=0pt,postaction={pattern=north east lines, pattern color=white}] (axis cs:3,2.06071428571429e-14) rectangle (axis cs:3.25,6.04821428571429e-14);
\draw[draw=white,fill=peru19615582,line width=0pt,postaction={pattern=north east lines, pattern color=white}] (axis cs:4,1.67571428571429e-14) rectangle (axis cs:4.25,4.88285714285714e-14);
\draw[draw=white,fill=peru19615582,line width=0pt,postaction={pattern=north east lines, pattern color=white}] (axis cs:5,1.3952380952381e-14) rectangle (axis cs:5.25,4.14166666666667e-14);
\draw[draw=white,fill=peru17813196,line width=0pt,postaction={pattern=north east lines, pattern color=white}] (axis cs:0,2.4e-13) rectangle (axis cs:0.25,3.15285714285714e-13);
\draw[draw=white,fill=peru17813196,line width=0pt,postaction={pattern=north east lines, pattern color=white}] (axis cs:1,1.21428571428571e-13) rectangle (axis cs:1.25,1.59928571428571e-13);
\draw[draw=white,fill=peru17813196,line width=0pt,postaction={pattern=north east lines, pattern color=white}] (axis cs:2,7.96428571428571e-14) rectangle (axis cs:2.25,1.05095238095238e-13);
\draw[draw=white,fill=peru17813196,line width=0pt,postaction={pattern=north east lines, pattern color=white}] (axis cs:3,6.04821428571429e-14) rectangle (axis cs:3.25,7.95357142857143e-14);
\draw[draw=white,fill=peru17813196,line width=0pt,postaction={pattern=north east lines, pattern color=white}] (axis cs:4,4.88285714285714e-14) rectangle (axis cs:4.25,6.41857142857143e-14);
\draw[draw=white,fill=peru17813196,line width=0pt,postaction={pattern=north east lines, pattern color=white}] (axis cs:5,4.14166666666667e-14) rectangle (axis cs:5.25,5.44285714285714e-14);
\draw[draw=white,fill=mediumpurple148107178,line width=0pt,postaction={pattern=north east lines, pattern color=white}] (axis cs:0,3.15285714285714e-13) rectangle (axis cs:0.25,1.62245714285714e-11);
\draw[draw=white,fill=mediumpurple148107178,line width=0pt,postaction={pattern=north east lines, pattern color=white}] (axis cs:1,1.59928571428571e-13) rectangle (axis cs:1.25,1.41560357142857e-11);
\draw[draw=white,fill=mediumpurple148107178,line width=0pt,postaction={pattern=north east lines, pattern color=white}] (axis cs:2,1.05095238095238e-13) rectangle (axis cs:2.25,1.41425e-11);
\draw[draw=white,fill=mediumpurple148107178,line width=0pt,postaction={pattern=north east lines, pattern color=white}] (axis cs:3,7.95357142857143e-14) rectangle (axis cs:3.25,1.4077625e-11);
\draw[draw=white,fill=mediumpurple148107178,line width=0pt,postaction={pattern=north east lines, pattern color=white}] (axis cs:4,6.41857142857143e-14) rectangle (axis cs:4.25,1.41229571428571e-11);
\draw[draw=white,fill=mediumpurple148107178,line width=0pt,postaction={pattern=north east lines, pattern color=white}] (axis cs:5,5.44285714285714e-14) rectangle (axis cs:5.25,1.40694285714286e-11);
\draw[draw=white,fill=indianred166106104,line width=0pt,postaction={pattern=north east lines, pattern color=white}] (axis cs:0,1.62245714285714e-11) rectangle (axis cs:0.25,1.97331428571429e-11);
\draw[draw=white,fill=indianred166106104,line width=0pt,postaction={pattern=north east lines, pattern color=white}] (axis cs:1,1.41560357142857e-11) rectangle (axis cs:1.25,1.781625e-11);
\draw[draw=white,fill=indianred166106104,line width=0pt,postaction={pattern=north east lines, pattern color=white}] (axis cs:2,1.41425e-11) rectangle (axis cs:2.25,1.79820238095238e-11);
\draw[draw=white,fill=indianred166106104,line width=0pt,postaction={pattern=north east lines, pattern color=white}] (axis cs:3,1.4077625e-11) rectangle (axis cs:3.25,1.80531785714286e-11);
\draw[draw=white,fill=indianred166106104,line width=0pt,postaction={pattern=north east lines, pattern color=white}] (axis cs:4,1.41229571428571e-11) rectangle (axis cs:4.25,1.81773857142857e-11);
\draw[draw=white,fill=indianred166106104,line width=0pt,postaction={pattern=north east lines, pattern color=white}] (axis cs:5,1.40694285714286e-11) rectangle (axis cs:5.25,1.81796309523809e-11);
\draw[draw=white,fill=gray152122143,line width=0pt,postaction={pattern=north east lines, pattern color=white}] (axis cs:0,1.97331428571429e-11) rectangle (axis cs:0.25,3.00412857142857e-11);
\draw[draw=white,fill=gray152122143,line width=0pt,postaction={pattern=north east lines, pattern color=white}] (axis cs:1,1.781625e-11) rectangle (axis cs:1.25,2.80204285714286e-11);
\draw[draw=white,fill=gray152122143,line width=0pt,postaction={pattern=north east lines, pattern color=white}] (axis cs:2,1.79820238095238e-11) rectangle (axis cs:2.25,2.81482380952381e-11);
\draw[draw=white,fill=gray152122143,line width=0pt,postaction={pattern=north east lines, pattern color=white}] (axis cs:3,1.80531785714286e-11) rectangle (axis cs:3.25,2.81879642857143e-11);
\draw[draw=white,fill=gray152122143,line width=0pt,postaction={pattern=north east lines, pattern color=white}] (axis cs:4,1.81773857142857e-11) rectangle (axis cs:4.25,2.82801428571429e-11);
\draw[draw=white,fill=gray152122143,line width=0pt,postaction={pattern=north east lines, pattern color=white}] (axis cs:5,1.8179630952381e-11) rectangle (axis cs:5.25,2.82729880952381e-11);
\draw[draw=white,fill=darkkhaki213187103,line width=0pt,postaction={pattern=north east lines, pattern color=white}] (axis cs:0,3.00412857142857e-11) rectangle (axis cs:0.25,3.17984285714286e-11);
\draw[draw=white,fill=darkkhaki213187103,line width=0pt,postaction={pattern=north east lines, pattern color=white}] (axis cs:1,2.80204285714286e-11) rectangle (axis cs:1.25,2.97357857142857e-11);
\draw[draw=white,fill=darkkhaki213187103,line width=0pt,postaction={pattern=north east lines, pattern color=white}] (axis cs:2,2.81482380952381e-11) rectangle (axis cs:2.25,2.98485476190476e-11);
\draw[draw=white,fill=darkkhaki213187103,line width=0pt,postaction={pattern=north east lines, pattern color=white}] (axis cs:3,2.81879642857143e-11) rectangle (axis cs:3.25,2.9879375e-11);
\draw[draw=white,fill=darkkhaki213187103,line width=0pt,postaction={pattern=north east lines, pattern color=white}] (axis cs:4,2.82801428571429e-11) rectangle (axis cs:4.25,2.99636285714286e-11);
\draw[draw=white,fill=darkkhaki213187103,line width=0pt,postaction={pattern=north east lines, pattern color=white}] (axis cs:5,2.82729880952381e-11) rectangle (axis cs:5.25,2.99538571428571e-11);
\draw[draw=white,fill=gray,very thin] (axis cs:-0.4,0) rectangle (axis cs:0.4,0);
\draw[draw=white,fill=gray,very thin,postaction={pattern=north east lines, pattern color=white}] (axis cs:-0.4,0) rectangle (axis cs:0.4,0);

  \end{axis}
\node [fill opacity=0.8,
    draw opacity=1,
    text opacity=1,
    draw=lightgray204,
    fill=gainsboro229,
    yshift=-1cm,
    anchor=west,name=leg2] at (leg.south west) {\shortstack[l]{
    \adjustbox{cframe=gray, bgcolor=gray}{\textcolor{gray}{\bf{\tiny{///}}}} $\blocktoep\,\,\,$ Matvec \\
\adjustbox{cframe=gray, bgcolor=gray}{\textcolor{white}{\bf{\tiny{///}}}} $\blocktoep^*$ Matvec }};

\end{tikzpicture}

%% file: images/new-single-gpu-scaling-line.tex
\begin{tikzpicture}

\definecolor{chocolate2267451}{RGB}{226,74,51}
\definecolor{darkslategray38}{RGB}{38,38,38}
\definecolor{gray119}{RGB}{119,119,119}
\definecolor{lavender234234242}{RGB}{234,234,242}
\definecolor{lightgray204}{RGB}{204,204,204}
\definecolor{mediumpurple152142213}{RGB}{152,142,213}
\definecolor{sandybrown25119394}{RGB}{251,193,94}
\definecolor{steelblue52138189}{RGB}{52,138,189}
  \definecolor{gainsboro229}{RGB}{229,229,229}

\definecolor{lightpink255181184}{RGB}{255,181,184}

\definecolor{yellowgreen14218666}{RGB}{142,186,66}

\begin{axis}[
title={Total Time ($\blocktoep$ Matvec)},
axis background/.style={fill=gainsboro229},
name=ax1,
axis line style={white},
legend cell align={left},
legend style={
  fill opacity=0.8,
  draw opacity=1,
  text opacity=1,
  at={(0.03,0.97)},
  anchor=north west,
  draw=lightgray204,
  fill=gainsboro229
},
tick align=outside,
tick pos=left,
scaled ticks=false,
minor y tick num=3,
y tick label style={/pgf/number format/fixed},
x grid style={white},
xlabel=\textcolor{darkslategray38}{Spatial parameter dimension $\numparam$ \(\displaystyle \qty(\times 10^3)\)},
xmajorgrids,
xmin=3.25, xmax=41.75,
xtick style={color=darkslategray38},
xtick={10,20,30,40},
xticklabels={1,2,3,4},
y grid style={white},
ylabel=\textcolor{darkslategray38}{Time $\qty(10^{-3}\text{s})$},
ymajorgrids,
ymajorticks=true,
ymin=-0.0001095205, ymax=0.0115566505,
ytick={0,0.002,0.004,0.006,0.008,0.01,0.012},
yticklabels={0,2,4,6,8,10,12},
ytick style={color=darkslategray38}
]
\addplot [semithick, steelblue52138189, mark=*, mark size=3, mark options={solid}]
table {%
5 0.000547528266906738
10 0.000907659530639648
15 0.00129413604736328
20 0.00163698196411133
25 0.00197041034698486
30 0.00233268737792969
35 0.00268685817718506
40 0.00301766395568848
};
\addlegendentry{20 Sensors}
\addplot [semithick, sandybrown25119394, mark=square*, mark size=3, mark options={solid}]
table {%
5 0.000754952430725098
10 0.00135695934295654
15 0.00196719169616699
20 0.00252532958984375
25 0.00312173366546631
30 0.00368154048919678
35 0.00428712368011475
40 0.00484812259674072
};
\addlegendentry{40 Sensors}
\addplot [semithick, yellowgreen14218666, mark=triangle*, mark size=3, mark options={solid}]
table {%
5 0.00101256370544434
10 0.00188577175140381
15 0.00267302989959717
20 0.0036168098449707
25 0.00442337989807129
30 0.0052412748336792
35 0.00620543956756592
40 0.00712895393371582
};
\addlegendentry{60 Sensors}
\addplot [semithick, chocolate2267451, mark=diamond*, mark size=3, mark options={solid}]
table {%
5 0.00117838382720947
10 0.00230884552001953
15 0.00336706638336182
20 0.00442469120025635
25 0.00555646419525146
30 0.00661730766296387
35 0.00771808624267578
40 0.00869810581207275
};
\addlegendentry{80 Sensors}
\addplot [semithick, mediumpurple152142213, mark=+, mark size=3, mark options={solid}]
table {%
5 0.00149333477020264
10 0.00284397602081299
15 0.00412189960479736
20 0.00551843643188477
25 0.00699329376220703
30 0.00832986831665039
35 0.00926387310028076
40 0.0106481313705444
};
\addlegendentry{100 Sensors}
\legend{};
\end{axis}

\begin{axis}[
  title={Time per Matrix Element},
  at={(ax1.south east)},
  name=ax2,
  xshift=1.75cm,
  axis background/.style={fill=gainsboro229},
  hide scale,
  axis line style={white},
  legend pos=outer north east,
  legend cell align={left},
  legend style={
    fill opacity=0.8,
    draw opacity=1,
    text opacity=1,
    draw=lightgray204,
    fill=gainsboro229
  },
  minor y tick num=3,
  tick align=outside,
  tick pos=left,
  x grid style={white},
  xlabel=\textcolor{darkslategray38}{Spatial parameter dimension $\numparam$ \(\displaystyle \qty(\times 10^3)\)},
  xmajorgrids,
  xmin=3.25, xmax=41.75,
xtick style={color=darkslategray38},
xtick={10,20,30,40},
xticklabels={1,2,3,4},
y grid style={white},
ylabel=\textcolor{darkslategray38}{Normalized Time $\qty(10^{-11}\text{s})$},
ymajorgrids,
scaled ticks=false,
ymajorticks=true,
ymin=1.179202e-11, ymax=3.2e-11,
ytick style={color=darkslategray38},
ytick={1e-11,1.5e-11,2e-11,2.5e-11,3e-11,3.5e-11,4e-11,4.5e-11},
yticklabels={1.0,1.5,2.0,2.5,3.0,3.5,4.0,4.5}
]
\addplot [semithick, steelblue52138189, mark=*, mark size=3, mark options={solid}]
table {%
5 2.73735e-11
10 2.2691e-11
15 2.15696666666667e-11
20 2.0463e-11
25 1.97036e-11
30 1.94390833333333e-11
35 1.91917142857143e-11
40 1.8860625e-11
};
\addlegendentry{20 Sensors}
\addplot [semithick, sandybrown25119394, mark=square*, mark size=3, mark options={solid}]
table {%
5 1.88735e-11
10 1.6962e-11
15 1.639375e-11
20 1.5783e-11
25 1.560875e-11
30 1.53398333333333e-11
35 1.53113214285714e-11
40 1.515034375e-11
};
\addlegendentry{40 Sensors}
\addplot [semithick, yellowgreen14218666, mark=triangle*, mark size=3, mark options={solid}]
table {%
5 1.68758333333333e-11
10 1.57146666666667e-11
15 1.48500555555556e-11
20 1.50699583333333e-11
25 1.47446333333333e-11
30 1.45591666666667e-11
35 1.47747380952381e-11
40 1.48519791666667e-11
};
\addlegendentry{60 Sensors}
\addplot [semithick, chocolate2267451, mark=diamond*, mark size=3, mark options={solid}]
table {%
5 1.47295e-11
10 1.44304375e-11
15 1.40292083333333e-11
20 1.382728125e-11
25 1.3891175e-11
30 1.3786e-11
35 1.37823928571429e-11
40 1.3590796875e-11
};
\addlegendentry{80 Sensors}
\addplot [semithick, mediumpurple152142213, mark=+, mark size=3, mark options={solid}]
table {%
5 1.49332e-11
10 1.422015e-11
15 1.37397333333333e-11
20 1.3796075e-11
25 1.398662e-11
30 1.38832e-11
35 1.32341142857143e-11
40 1.33101e-11
};
\addlegendentry{100 Sensors}
  \end{axis}
  
\end{tikzpicture}

%% file: images/new-single-gpu-time-scaling-line.tex
\begin{tikzpicture}

\definecolor{chocolate2267451}{RGB}{226,74,51}
\definecolor{darkslategray38}{RGB}{38,38,38}
\definecolor{lavender234234242}{RGB}{234,234,242}
\definecolor{lightgray204}{RGB}{204,204,204}
  \definecolor{gainsboro229}{RGB}{229,229,229}
\definecolor{steelblue52138189}{RGB}{52,138,189}

\begin{axis}[
title={Total Time},
axis background/.style={fill=gainsboro229},
name=ax1,
axis line style={white},
legend cell align={left},
legend style={
  fill opacity=0.8,
  draw opacity=1,
  text opacity=1,
  at={(0.03,0.97)},
  anchor=north west,
  draw=lightgray204,
  fill=gainsboro229
},
y tick label style={/pgf/number format/fixed},
tick align=outside,
x grid style={white},
xlabel=\textcolor{darkslategray38}{Time Steps \(\displaystyle \numtime \ \qty(\times 10^3)\)},
xmajorgrids,
xmajorticks,
tick align=outside,
minor y tick num=3,
tick pos=left,
xmin=5.5, xmax=104.5,
xtick style={color=darkslategray38},
xtick={10,20,30,40,50,60,70,80,90,100},
xticklabels={1,2,3,4,5,6,7,8,9,10},
scaled ticks=false,
y grid style={white},
ylabel=\textcolor{darkslategray38}{Time $\qty(10^{-3}\text{s})$},
ymajorgrids,
ymajorticks=true,
ymin=0.0006328955, ymax=0.0112174145,
ytick style={color=darkslategray38},
ytick={0,0.002,0.004,0.006,0.008,0.01,0.012},
yticklabels={0,2,4,6,8,10,12}
]
\addplot [semithick, chocolate2267451, mark=square*, mark size=3, mark options={solid}]
table {%
10 0.00118006
20 0.00224235
30 0.00340936
40 0.00439925
50 0.00546915
60 0.00665626
70 0.0077179
80 0.00852665
90 0.00957035
100 0.01043149
};
\addlegendentry{$\blocktoep$ Matvec}
\addplot [semithick, steelblue52138189, mark=*, mark size=3, mark options={solid}]
table {%
10 0.00111401
20 0.00217895
30 0.00336757
40 0.00438165
50 0.00551234
60 0.00667784
70 0.00776856
80 0.00865137
90 0.00983098
100 0.0107363
};
\addlegendentry{$\blocktoep^*$ Matvec}
\end{axis}

\begin{axis}[
  title={Time per Matrix Element},
  at={(ax1.south east)},
  hide scale,
  name=ax2,
  xshift=2cm,
  axis background/.style={fill=gainsboro229},
  axis line style={white},
  legend cell align={left},
  legend style={
    fill opacity=0.8,
    draw opacity=1,
    text opacity=1,
    draw=lightgray204,
    fill=gainsboro229,
    at={(0.97,0.97)},
    anchor=north east,
  },
  tick align=outside,
  x grid style={white},
  xlabel=\textcolor{darkslategray38}{Time Steps \(\displaystyle \numtime\ \qty(\times 10^3)\) },
  xmajorgrids,
  xmajorticks,
  tick align=outside,
  minor y tick num=3,
  tick pos=left,
  xmin=5.5, xmax=104.5,
xtick style={color=darkslategray38},
xtick={10,20,30,40,50,60,70,80,90,100},
xticklabels={1,2,3,4,5,6,7,8,9,10},
  xtick style={color=darkslategray38},
  y grid style={white},
  ylabel=\textcolor{darkslategray38}{Time $\qty(10^{-11}\text{s})$},
  ymajorgrids,
  ymajorticks,
  ymin=1.2953793125e-11, ymax=1.4836319375e-11,
ytick style={color=darkslategray38},
ytick={1.275e-11,1.3e-11,1.325e-11,1.35e-11,1.375e-11,1.4e-11,1.425e-11,1.45e-11,1.475e-11,1.5e-11},
yticklabels={1.275,1.300,1.325,1.350,1.375,1.400,1.425,1.450,1.475,1.500}
]
\addplot [semithick, chocolate2267451, mark=square*, mark size=3, mark options={solid}]
table {%
10 1.475075e-11
20 1.40146875e-11
30 1.42056666666667e-11
40 1.374765625e-11
50 1.3672875e-11
60 1.38672083333333e-11
70 1.37819642857143e-11
80 1.3322890625e-11
90 1.32921527777778e-11
100 1.30393625e-11
};
\addlegendentry{$\blocktoep$ Matvec}
\addplot [semithick, steelblue52138189, mark=*, mark size=3, mark options={solid}]
table {%
10 1.3925125e-11
20 1.36184375e-11
30 1.40315416666667e-11
40 1.369265625e-11
50 1.378085e-11
60 1.39121666666667e-11
70 1.38724285714286e-11
80 1.3517765625e-11
90 1.36541388888889e-11
100 1.3420375e-11
};
\addlegendentry{$\blocktoep^*$ Matvec}
\end{axis}

\end{tikzpicture}

%% file: images/new_roofline.tex
\begin{tikzpicture}

\definecolor{crimson2143940}{RGB}{214,39,40}
\definecolor{darkorange25512714}{RGB}{255,127,14}
\definecolor{darkturquoise23190207}{RGB}{23,190,207}
\definecolor{dimgray85}{RGB}{85,85,85}
\definecolor{forestgreen4416044}{RGB}{44,160,44}
\definecolor{gainsboro229}{RGB}{229,229,229}
\definecolor{goldenrod18818934}{RGB}{188,189,34}
\definecolor{lightgray204}{RGB}{204,204,204}
\definecolor{mediumpurple148103189}{RGB}{148,103,189}
\definecolor{orchid227119194}{RGB}{227,119,194}
\definecolor{sienna1408675}{RGB}{140,86,75}
\definecolor{steelblue31119180}{RGB}{31,119,180}
\definecolor{yellow}{RGB}{255,255,0}

\begin{axis}[
axis background/.style={fill=gainsboro229},
axis line style={white},
title={Roofline Analysis},
legend cell align={left},
legend style={
  fill opacity=0.8,
  draw opacity=1,
  text opacity=1,
  at={(1.01,0)},
  anchor=south west,
  draw=lightgray204,
  fill=gainsboro229
},
legend pos=outer north east,
log basis x={10},
log basis y={10},
tick align=outside,
tick pos=left,
x grid style={white},
xlabel={Arithmetic intensity (FLOPs/byte)},
xmajorgrids,
xmin=0.00707945784384138, xmax=14.1253754462275,
xmode=log,
xtick style={color=dimgray85},
xtick={0.0001,0.001,0.01,0.1,1,10,100,1000},
xticklabels={
  \(\displaystyle {10^{-4}}\),
  \(\displaystyle {10^{-3}}\),
  \(\displaystyle {10^{-2}}\),
  \(\displaystyle {10^{-1}}\),
  \(\displaystyle {10^{0}}\),
  \(\displaystyle {10^{1}}\),
  \(\displaystyle {10^{2}}\),
  \(\displaystyle {10^{3}}\)
},
y grid style={white},
ylabel={Performance (FLOPs/s)},
ymajorgrids,
ymin=1437269267.294, ymax=9869621001549.32,
ymode=log,
ytick style={color=dimgray85},
ytick={1000000000,10000000000,100000000000,1000000000000,10000000000000,100000000000000},
yticklabels={
  \(\displaystyle {10^{9}}\),
  \(\displaystyle {10^{10}}\),
  \(\displaystyle {10^{11}}\),
  \(\displaystyle {10^{12}}\),
  \(\displaystyle {10^{13}}\),
  \(\displaystyle {10^{14}}\)
}
]
\addplot [semithick, black]
table {%
0.01 19341131690.7
3.79 7334267289273.66
10 7334267289273.66
};
\addlegendentry{Roofline}
\addplot [semithick, steelblue31119180, mark=square*, mark size=3, mark options={solid}, only marks]
table {%
0.13 205131185383.73
};
\addlegendentry{$\blocktoep$: FFT$\qty(\paramvec)$}
\addplot [semithick, darkorange25512714, mark=square*, mark size=3, mark options={solid,rotate=180}, only marks]
table {%
0.5 786822150464.23
};
\addlegendentry{$\blocktoep$: SBGEMV}
\addplot [semithick, mediumpurple148103189, mark=*, mark size=3, mark options={solid}, only marks]
table {%
0.26 10322580645.16
};
\addlegendentry{$\blocktoep$: IFFT$\qty(\datavec)$}
\addplot [semithick, orchid227119194, mark=triangle*, mark size=3, mark options={solid}, only marks]
table {%
0.26 10191082802.55
};
\addlegendentry{$\blocktoep^*$: FFT$\qty(\datavec)$}
\addplot [semithick, sienna1408675, mark=triangle*, mark size=3, mark options={solid}, only marks]
table {%
0.51 834978912211.41
};
\addlegendentry{$\blocktoep^*$: SBGEMV}
\addplot [semithick, goldenrod18818934, mark=diamond*, mark size=3, mark options={solid}, only marks]
table {%
0.13 203629969113.41
};
\addlegendentry{$\blocktoep^*$: IFFT$\qty(\paramvec)$}
\end{axis}

\end{tikzpicture}

%% file: images/new_weak_scaling_full_ls6.tex
\begin{tikzpicture}

\definecolor{chocolate2267451}{RGB}{226,74,51}
\definecolor{dimgray85}{RGB}{85,85,85}
\definecolor{gainsboro229}{RGB}{229,229,229}
\definecolor{lightgray204}{RGB}{204,204,204}
\definecolor{mediumpurple152142213}{RGB}{152,142,213}
\definecolor{steelblue52138189}{RGB}{52,138,189}
\definecolor{chocolate2267451}{RGB}{226,74,51}
\definecolor{darkslategray38}{RGB}{38,38,38}
\definecolor{lavender234234242}{RGB}{234,234,242}
\definecolor{lightgray204}{RGB}{204,204,204}
\definecolor{mediumpurple152142213}{RGB}{152,142,213}
\definecolor{steelblue52138189}{RGB}{52,138,189}

\begin{axis}[
title={Strong Scaling},
axis background/.style={fill=gainsboro229},
axis line style={white},
name=ax1,
legend cell align={left},
legend style={
  fill opacity=0.8,
  draw opacity=1,
  text opacity=1,
  at={(0.03,0.97)},
  anchor=north west,
  draw=lightgray204,
  fill=gainsboro229
},
tick align=outside,
tick pos=left,
minor y tick num=3,
x grid style={white},
xlabel={Number of GPUs},
xmajorgrids,
xmin=5.5, xmax=49.5,
xtick style={color=dimgray85},
minor x tick num=5,
xtick={6,12,18,24,30,36,42,48},
y grid style={white},
ylabel={Speedup},
ymajorgrids,
ymin=0.75, ymax=8.25,
ytick style={color=dimgray85}
]
\addplot [semithick, chocolate2267451, mark=square*, mark size=3, mark options={solid}]
table {%
6 1
12 1.92917609214783
18 2.72568559646606
24 3.35736203193665
30 3.7843451499939
36 4.09403228759766
42 4.3079309463501
48 4.44393968582153
};
\addlegendentry{Real Speedup $\blocktoep$}
\addplot [semithick, steelblue52138189, mark=*, mark size=3, mark options={solid}]
table {%
6 1
12 1.91036021709442
18 2.66540145874023
24 3.25368523597717
30 3.6377124786377
36 3.85931777954102
42 4.02374887466431
48 4.09837913513184
};
\addlegendentry{Real Speedup $\blocktoep^*$}
\addplot [thick, dashed, mediumpurple152142213]
table {%
6 1
48 8
};
\addlegendentry{Ideal Speedup}
\end{axis}

\hskip 5pt

\begin{axis}[
  title={Weak Scaling},
  at={(ax1.south east)},
  name=ax2,
  xshift=1.8cm,
  axis background/.style={fill=gainsboro229},
  axis line style={white},
  legend cell align={left},
  legend style={
    fill opacity=0.8,
    draw opacity=1,
    text opacity=1,
    at={(0.03,0.03)},
    anchor=south west,
    draw=lightgray204,
    fill=gainsboro229
  },
  minor y tick num=3,
  tick align=outside,
  tick pos=left,
  x grid style={white},
xlabel=\textcolor{darkslategray38}{Number of GPUs},
xmajorgrids,
xmin=3.9, xmax=50.1,
minor x tick num=5,
xtick={6,12,18,24,30,36,42,48},
xtick style={color=darkslategray38},
y grid style={white},
ylabel=\textcolor{darkslategray38}{Efficiency},
ymajorgrids,
ymin=0, ymax=1.05,
ytick style={color=darkslategray38}
  ]
  \addplot [semithick, chocolate2267451, mark=square*, mark size=3, mark options={solid}]
table {%
6 1
12 0.977975726127625
18 0.95905065536499
24 0.940670728683472
30 0.92168402671814
36 0.905966758728027
42 0.889054298400879
48 0.872745990753174
};
\addlegendentry{Weak Scaling Efficiency $\blocktoep$}
\addplot [semithick, steelblue52138189, mark=*, mark size=3, mark options={solid}]
table {%
6 1
12 0.976243734359741
18 0.954857110977173
24 0.934655547142029
30 0.912842035293579
36 0.894962191581726
42 0.876702666282654
48 0.857527732849121
};
\addlegendentry{Weak Scaling Efficiency $\blocktoep^*$}
\addplot [thick, dashed, mediumpurple152142213]
table {%
6 1
48 1
};
\addlegendentry{Ideal Efficiency}
  \end{axis}

\end{tikzpicture}

%% file: images/grid_test.tex
\begin{tikzpicture}

\definecolor{chocolate2267451}{RGB}{226,74,51}
\definecolor{dimgray85}{RGB}{85,85,85}
\definecolor{gainsboro229}{RGB}{229,229,229}
\definecolor{lightgray204}{RGB}{204,204,204}
\definecolor{mediumpurple152142213}{RGB}{152,142,213}
\definecolor{steelblue52138189}{RGB}{52,138,189}

\begin{axis}[
axis background/.style={fill=gainsboro229},
axis line style={white},
legend cell align={left},
legend style={
  fill opacity=0.8,
  draw opacity=1,
  text opacity=1,
  at={(0.03,0.97)},
  anchor=north west,
  draw=lightgray204,
  fill=gainsboro229
},
tick align=outside,
tick pos=left,
x grid style={white},
xlabel={Number of Rows $\numprocrows$ in Processor Grid},
xmajorgrids,
xmin=0.65, xmax=8.35,
xtick={1,2,4,6,8},
xtick style={color=dimgray85},
y grid style={white},
minor y tick num=3,
ylabel={Communication Time (normalized)},
ymajorgrids,
ymin=0.768607554879182, ymax=5.85924134753718,
ytick style={color=dimgray85}
]
\addplot [semithick, dashed, chocolate2267451, mark=square*, mark size=3, mark options={solid}]
table {%
1 1
2 1.26372897624969
4 2.12872123718262
5 2.61080527305603
8 5.62784910202026
};
\addlegendentry{$\log \frac{\numdata}{\numparam}=-4$}
\addplot [semithick, dashed, steelblue52138189, mark=*, mark size=3, mark options={solid}]
table {%
1 1
2 1.00444221496582
4 1.36001825332642
5 1.63141071796417
8 3.2625937461853
};
\addlegendentry{$\log \frac{\numdata}{\numparam}=-3$}
\addplot [semithick, dashed, mediumpurple152142213, mark=diamond*, mark size=3, mark options={solid}]
table {%
1 1.28495573997498
2 1.120232462883
4 1
5 1.05418395996094
8 1.59659624099731
};
\addlegendentry{$\log \frac{\numdata}{\numparam}=-2$}
\end{axis}

\end{tikzpicture}

%% file: conclusions.tex
\revstart
In the context of solution of linear inverse problems for high-dimensional parameter fields, the Hessian action on a vector typically has to be performed a large number of times. 
\revend
This is especially true for Hessian matrices with a slow spectral decay, i.e.~where the effective rank is large compared to the parameter or data dimension, which prevents the effective use of low-rank based techniques. In such cases, and particularly for real-time inversion, it is essential that the Hessian matvec be performed both rapidly and efficiently.

In this paper, we proposed efficient and scalable multi-GPU
FFT-accelerated algorithms for performing matvecs of block Toeplitz
systems that arise in inverse problems involving discrete
shift-invariant systems. The numerical results illustrate the
applicability of our approach to Hessian matvecs for large-scale
inverse problems governed by autonomous systems. For time-invariant
dynamical systems, the cost of the Hessian matvec increases only
linearly with the number of time steps, and near-linearly in its
one-time setup cost, instead of quadratically for a naive
implementation. While classical inverse algorithms with adjoint-based
Hessian matvecs also scale linearly with the number of time steps, the
overall cost per Hessian matvec is several orders of magnitude larger
for many problems, since each Hessian matvec requires one forward and
one adjoint PDE solve.

This cost becomes even more severe when (1)~the governing PDEs are
discretized with high-order methods, (2)~the time-stepping methods for
the PDE solves require small time steps for stability or accuracy
relative to the data acquisition rate, or (3)~the governing PDE system describes a multiphysics problem that is discretized with many DOFs per grid point. In those cases, our approach can outperform the classical algorithm by a factor of $1{,}000\times$ or larger, depending on the particular problem. Even for a problem that did not meet the prior three criteria, we showed that our approach was over $750\times$ more efficient (measured in FLOPs).

Our algorithm runs efficiently on GPUs, achieving more than 80 percent
of peak bandwidth on an NVIDIA A100 GPU. Parallel scalability of the
algorithm is limited only by the scaling of a vector broadcast and
reduction operation that is needed once per matvec, achieving 85-87\% weak scaling efficiency on 48 A100 GPUs.

\revstart 
The GPU-accelerated FFT-based block Toeplitz matvec algorithm introduced in this paper has been applied to a large-scale 3D implementation of a recently developed inverse model for tsunami early warning~\cite{henneking2025bell}. The algorithm proposed here is an essential building block for solving the tsunami inverse problem, which is governed by an acoustic--gravity wave propagation model representing an autonomous dynamical system~\cite{lotto2015tsunami}. This tsunami inversion was applied to real-time inference of seafloor motion and forecasting of tsunami wave heights for realistic rupture events in the Cascadia subduction zone.
\revend

%% file: appendix.tex
As mentioned in~\cref{sec:ReadingBlockToeplitz}, in SOTI ordering, $\widehat{\toepblock}$ is a matrix with $\numparam \times \numdata$ diagonal blocks. A natural algorithm to multiply this matrix with a vector (also in SOTI ordering) would be to compute elementwise products between corresponding matrix and vector blocks and then sum over block rows. This was in fact the original implementation of the algorithm. Elementwise products and reductions over local blocks of the matrix were computed with optimized custom CUDA kernels. To handle matvecs with $\blocktoep^*$, the only necessary changes were to take the complex conjugate of the matrix elements before computing elementwise products and to reduce over columns instead of rows.

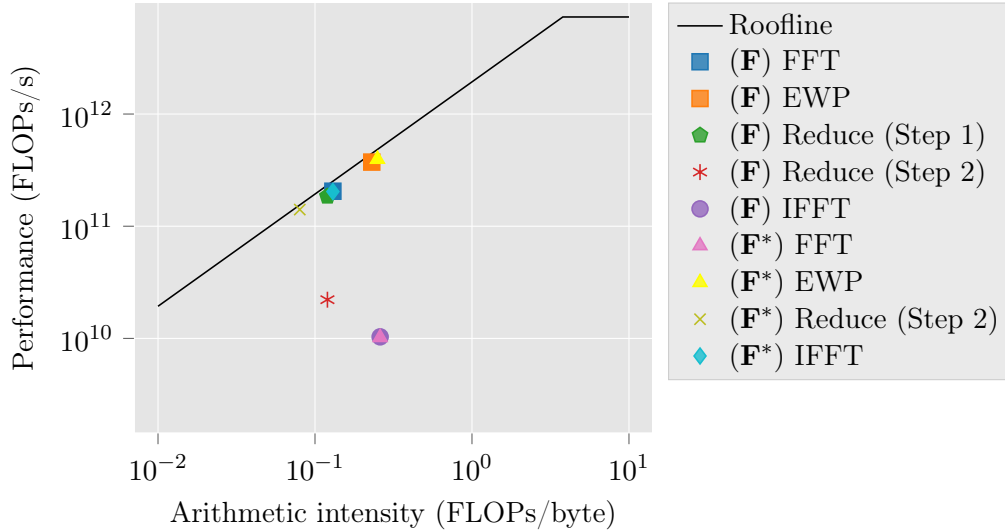
\begin{figure}[htpb]
    \centering
    \input{images/roofline.tex}
    \caption{Roofline plots for the main kernels used in the matvec computation (alternate algorithm). All major kernels (including custom kernels) are operating at close to peak performance. Some kernels are not operating at peak performance due to small input size (they operate on a vector of size $\procrowblock \ll \proccolblock$); these kernels comprise less than 1\% of the total runtime.}\label{fig:OldRoofline}
\end{figure}

\Cref{fig:OldRoofline} shows the roofline plot for this implementation of the matvec algorithm. All major kernels, including custom kernels, are operating at close to peak performance. Similarly to the plot in~\cref{fig:Roofline}, some kernels are operating under the roofline purely due to the small input size.

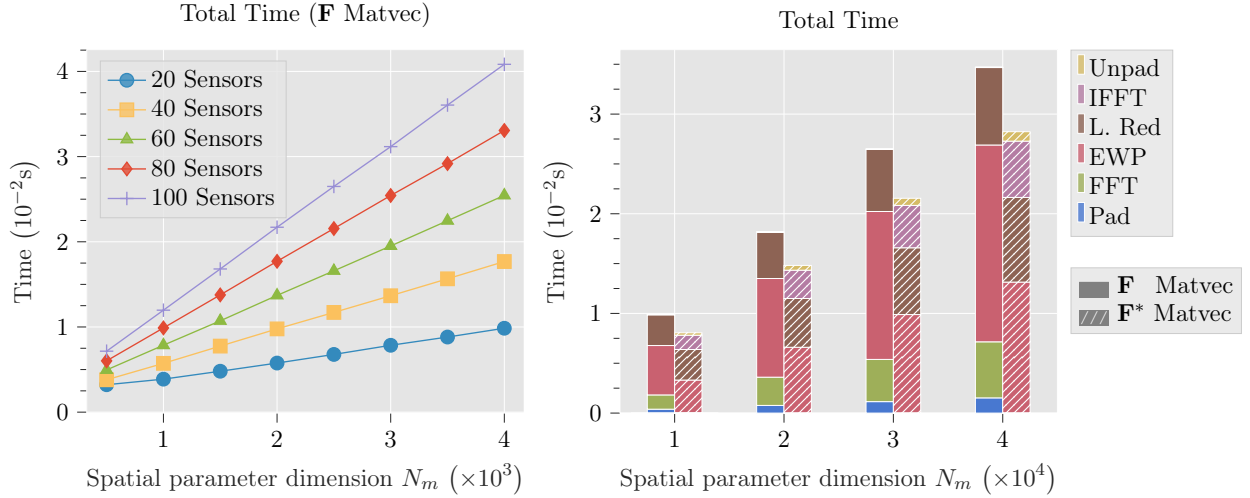
\begin{figure}[htbp]
  \centering
  \resizebox{\columnwidth}{!}{\input{images/single-gpu-line-mod}}
    \caption{Single GPU scaling results for the matvec (alternate algorithm). Here, $\numtime = 2{,}000$ for both plots. For the parameter scaling plot (right), $\numdata = 7$. The extra memory usage of the alternate algorithm means that only up to $\numparam = 40{,}000$ can be fit on the GPU. Compare to~\cref{fig:SingleGPUScaling,fig:SingleGPUScalingLine} (left). }\label{fig:OldSingleGPU}
\end{figure}

\Cref{fig:OldSingleGPU} depicts single GPU scaling results corresponding to~\cref{fig:SingleGPUScaling,fig:SingleGPUScalingLine} for the alternate algorithm. This figure shows that the local matvec computations are dominated by elementwise products (EWP) and local reductions (L. Red). Moreover, the overall algorithm runtimes are significantly longer than those of the algorithm presented in the main paper. The reason for this is the following: memory access patterns for the EWP are contiguous while those for the local reduction are strided. Strided access issues can be mitigated (and were done so) using grid-strided loops, as is evident from the fact that all kernels achieve 85--90\% memory bandwidth. However, the different memory access patterns for the EWP and local reduction imply that they cannot be effectively combined into the same kernel. As a result, data has to be transferred back and forth between device global memory and cache. For the case of the EWP and reduction, there are 3 matrix read/writes and 2 vector read/writes for the EWP$+$local reduction. In contrast,~\cref{alg:matvec,alg:matvec2} incur 1 matrix read/write and 4 vector read/writes. The strided access pattern is accounted for by the \texttt{SOTI\_TO\_TOSI} and \texttt{TOSI\_TO\_SOTI} reindexing operations (these also explain the additional vector reads). Overall, the reduced amount of matrix-sized I/O combined with the fact that the main computation can be offloaded to cuBLAS routines explains the superior performance of the algorithms in the main paper. In addition, the reindexing kernels operate only on vector-sized data; their runtime is negligible even though they attain a slightly lower percentage of peak bandwidth. Also, note that the alternate algorithm uses more memory overall for each matvec since intermediate results of the EWP have to be stored. As a result, the largest problem size that can fit on a single GPU is reduced for the alternate algorithm. This is why in~\cref{fig:OldSingleGPU} (right), only up to $\numparam=40{,}000$ is tested. 

%% file: images/roofline.tex
\begin{tikzpicture}

\definecolor{crimson2143940}{RGB}{214,39,40}
\definecolor{darkorange25512714}{RGB}{255,127,14}
\definecolor{darkturquoise23190207}{RGB}{23,190,207}
\definecolor{dimgray85}{RGB}{85,85,85}
\definecolor{forestgreen4416044}{RGB}{44,160,44}
\definecolor{gainsboro229}{RGB}{229,229,229}
\definecolor{goldenrod18818934}{RGB}{188,189,34}
\definecolor{lightgray204}{RGB}{204,204,204}
\definecolor{mediumpurple148103189}{RGB}{148,103,189}
\definecolor{orchid227119194}{RGB}{227,119,194}
\definecolor{sienna1408675}{RGB}{140,86,75}
\definecolor{steelblue31119180}{RGB}{31,119,180}
\definecolor{yellow}{RGB}{255,255,0}

\begin{axis}[
axis background/.style={fill=gainsboro229},
axis line style={white},
legend cell align={left},
legend style={
  fill opacity=0.8,
  draw opacity=1,
  text opacity=1,
  at={(1.01,0)},
  anchor=south west,
  draw=lightgray204,
  fill=gainsboro229
},
legend pos=outer north east,
log basis x={10},
log basis y={10},
tick align=outside,
tick pos=left,
x grid style={white},
xlabel={Arithmetic intensity (FLOPs/byte)},
xmajorgrids,
xmin=0.00707945784384138, xmax=14.1253754462275,
xmode=log,
xtick style={color=dimgray85},
xtick={0.0001,0.001,0.01,0.1,1,10,100,1000},
xticklabels={
  \(\displaystyle {10^{-4}}\),
  \(\displaystyle {10^{-3}}\),
  \(\displaystyle {10^{-2}}\),
  \(\displaystyle {10^{-1}}\),
  \(\displaystyle {10^{0}}\),
  \(\displaystyle {10^{1}}\),
  \(\displaystyle {10^{2}}\),
  \(\displaystyle {10^{3}}\)
},
y grid style={white},
ylabel={Performance (FLOPs/s)},
ymajorgrids,
ymin=1437269267.294, ymax=9869621001549.32,
ymode=log,
ytick style={color=dimgray85},
ytick={1000000000,10000000000,100000000000,1000000000000,10000000000000,100000000000000},
yticklabels={
  \(\displaystyle {10^{9}}\),
  \(\displaystyle {10^{10}}\),
  \(\displaystyle {10^{11}}\),
  \(\displaystyle {10^{12}}\),
  \(\displaystyle {10^{13}}\),
  \(\displaystyle {10^{14}}\)
}
]
\addplot [semithick, black]
table {%
0.01 19341131690.7
3.79 7334267289273.66
10 7334267289273.66
};
\addlegendentry{Roofline}
\addplot [semithick, steelblue31119180, mark=square*, mark size=3, mark options={solid}, only marks]
table {%
0.13 205131185383.73
};
\addlegendentry{($\blocktoep$) FFT}
\addplot [semithick, darkorange25512714, mark=square*, mark size=3, mark options={solid,rotate=180}, only marks]
table {%
0.23 372965911621.54
};
\addlegendentry{($\blocktoep$) EWP}
\addplot [semithick, forestgreen4416044, mark=pentagon*, mark size=3, mark options={solid}, only marks]
table {%
0.12 184551603581.19
};
\addlegendentry{($\blocktoep$) Reduce (Step 1)}
\addplot [semithick, crimson2143940, mark=asterisk, mark size=3, mark options={solid}, only marks]
table {%
0.12 22130270092.23
};
\addlegendentry{($\blocktoep$) Reduce (Step 2)}
\addplot [semithick, mediumpurple148103189, mark=*, mark size=3, mark options={solid}, only marks]
table {%
0.26 10322580645.16
};
\addlegendentry{($\blocktoep$) IFFT}
\addplot [semithick, orchid227119194, mark=triangle*, mark size=3, mark options={solid}, only marks]
table {%
0.26 10191082802.55
};
\addlegendentry{($\blocktoep^*$) FFT}
\addplot [semithick, yellow, mark=triangle*, mark size=3, mark options={solid}, only marks]
table {%
0.25 395332515585.68
};
\addlegendentry{($\blocktoep^*$) EWP}
\addplot [semithick, goldenrod18818934, mark=x, mark size=3, mark options={solid}, only marks]
table {%
0.08 140709040035.41
};
\addlegendentry{($\blocktoep^*$) Reduce (Step 2)}
\addplot [semithick, darkturquoise23190207, mark=diamond*, mark size=3, mark options={solid}, only marks]
table {%
0.13 203629969113.41
};
\addlegendentry{($\blocktoep^*$) IFFT}
\end{axis}

\end{tikzpicture}

%% file: images/single-gpu-line-mod.tex
\begin{tikzpicture}

\definecolor{chocolate2267451}{RGB}{226,74,51}
\definecolor{darkslategray38}{RGB}{38,38,38}
\definecolor{gray119}{RGB}{119,119,119}
\definecolor{lavender234234242}{RGB}{234,234,242}
\definecolor{lightgray204}{RGB}{204,204,204}
\definecolor{mediumpurple152142213}{RGB}{152,142,213}
\definecolor{sandybrown25119394}{RGB}{251,193,94}
\definecolor{steelblue52138189}{RGB}{52,138,189}
  \definecolor{gainsboro229}{RGB}{229,229,229}

\definecolor{lightpink255181184}{RGB}{255,181,184}

\definecolor{yellowgreen14218666}{RGB}{142,186,66}

\definecolor{darkkhaki15817589}{RGB}{158,175,89}
\definecolor{darkkhaki213187103}{RGB}{213,187,103}
\definecolor{dimgray1419984}{RGB}{141,99,84}
\definecolor{gray}{RGB}{128,128,128}
\definecolor{indianred20197112}{RGB}{201,97,112}
\definecolor{rosybrown180124163}{RGB}{180,124,163}
\definecolor{royalblue72120208}{RGB}{72,120,208}
\definecolor{peru17813196}{RGB}{178,131,96}
\definecolor{peru19615582}{RGB}{196,155,82}
\definecolor{darkkhaki213187103}{RGB}{213,187,103}
\definecolor{gray}{RGB}{128,128,128}
\definecolor{gray152122143}{RGB}{152,122,143}
\definecolor{indianred166106104}{RGB}{166,106,104}
\definecolor{mediumpurple148107178}{RGB}{148,107,178}

\begin{axis}[
title={Total Time ($\blocktoep$ Matvec)},
axis background/.style={fill=gainsboro229},
name=ax1,
axis line style={white},
legend cell align={left},
legend style={
  fill opacity=0.8,
  draw opacity=1,
  text opacity=1,
  at={(0.03,0.97)},
  anchor=north west,
  draw=lightgray204,
  fill=gainsboro229
},
tick align=outside,
tick pos=left,
scaled ticks=false,
minor y tick num=3,
y tick label style={/pgf/number format/fixed},
x grid style={white},
xlabel=\textcolor{darkslategray38}{Spatial parameter dimension $\numparam$ \(\displaystyle \qty(\times 10^3)\)},
xmajorgrids,
xmin=3.25, xmax=41.75,
xtick style={color=darkslategray38},
xtick={10,20,30,40},
xticklabels={1,2,3,4},
y grid style={white},
ylabel=\textcolor{darkslategray38}{Time $\qty(10^{-2}\text{s})$},
ymajorgrids,
ymajorticks=true,
ymin=-0.0001095205, ymax=0.0425566505,
ytick={0,0.01,0.02,0.03,0.04,0.05},
yticklabels={0,1,2,3,4,5},
ytick style={color=darkslategray38}
]
\addplot [semithick, steelblue52138189, mark=*, mark size=3, mark options={solid}]
table {%
5 0.00321888
10 0.00387434
15 0.00480368
20 0.00576294
25 0.00677603
30 0.00784963
35 0.00881773
40 0.00984988
};
\addlegendentry{20 Sensors}
\addplot [semithick, sandybrown25119394, mark=square*, mark size=3, mark options={solid}]
table {%
5 0.00379913
10 0.00572861
15 0.00775435
20 0.0097743
25 0.01170709
30 0.01366603
35 0.01565396
40 0.01768134
};
\addlegendentry{40 Sensors}
\addplot [semithick, yellowgreen14218666, mark=triangle*, mark size=3, mark options={solid}]
table {%
5 0.00495795
10 0.00785931
15 0.01071977
20 0.0137054
25 0.01658946
30 0.01951071
35 0.0224731
40 0.02544502
};
\addlegendentry{60 Sensors}
\addplot [semithick, chocolate2267451, mark=diamond*, mark size=3, mark options={solid}]
table {%
5 0.00602229
10 0.00988602
15 0.01376263
20 0.01770734
25 0.02155279
30 0.02543309
35 0.02917804
40 0.03306829
};
\addlegendentry{80 Sensors}
\addplot [semithick, mediumpurple152142213, mark=+, mark size=3, mark options={solid}]
table {%
5 0.00715676
10 0.01197084
15 0.01680897
20 0.02170712
25 0.02649564
30 0.03117082
35 0.03604998
40 0.04083316
};
\addlegendentry{100 Sensors}
\end{axis}

\begin{axis}[
title={Total Time},
axis background/.style={fill=gainsboro229},
name=ax2,
at={(ax1.south east)},
  name=ax2,
  xshift=1.5cm,
axis line style={white},
legend cell align={left},
legend style={
  fill opacity=0.8,
  draw opacity=1,
  text opacity=1,
  at={(0.03,0.97)},
  anchor=north west,
  draw=lightgray204,
  fill=gainsboro229,
  name=leg
},
legend pos=outer north east,
reverse legend,
minor y tick num=3,
y tick label style={/pgf/number format/fixed},
tick align=outside,
tick pos=left,
x grid style={white},
xlabel=\textcolor{darkslategray38}{Spatial parameter dimension $\numparam$ \(\displaystyle \qty(\times 10^4)\)},
xmajorgrids,
xmin=-0.5, xmax=3.5,
xtick style={color=darkslategray38},
xtick={0,1,2,3},
xticklabels={1,2,3,4},
y grid style={white},
ylabel=\textcolor{darkslategray38}{Time $\qty(10^{-2}\text{s})$},
scaled ticks=false,
ymajorgrids,
ymin=0, ymax=0.036467365,
ytick={0,0.01,0.02,0.03,0.04},
yticklabels={0,1,2,3,4},
ytick style={color=darkslategray38}
]
\draw[draw=white,fill=royalblue72120208,line width=0pt] (axis cs:-0.25,0) rectangle (axis cs:0,0.00039092);
\addlegendimage{ybar,ybar legend,draw=white,fill=royalblue72120208,line width=0pt}
\addlegendentry{Pad}

\draw[draw=white,fill=royalblue72120208,line width=0pt] (axis cs:0.75,0) rectangle (axis cs:1,0.00076949);
\draw[draw=white,fill=royalblue72120208,line width=0pt] (axis cs:1.75,0) rectangle (axis cs:2,0.00114933);
\draw[draw=white,fill=royalblue72120208,line width=0pt] (axis cs:2.75,0) rectangle (axis cs:3,0.00152556);
\draw[draw=white,fill=darkkhaki15817589,line width=0pt] (axis cs:-0.25,0.00039092) rectangle (axis cs:0,0.00182246);
\addlegendimage{ybar,ybar legend,draw=white,fill=darkkhaki15817589,line width=0pt}
\addlegendentry{FFT}

\draw[draw=white,fill=darkkhaki15817589,line width=0pt] (axis cs:0.75,0.00076949) rectangle (axis cs:1,0.00360291);
\draw[draw=white,fill=darkkhaki15817589,line width=0pt] (axis cs:1.75,0.00114933) rectangle (axis cs:2,0.00538552);
\draw[draw=white,fill=darkkhaki15817589,line width=0pt] (axis cs:2.75,0.00152556) rectangle (axis cs:3,0.0071511);
\draw[draw=white,fill=indianred20197112,line width=0pt] (axis cs:-0.25,0.00182246) rectangle (axis cs:0,0.00677816);
\addlegendimage{ybar,ybar legend,draw=white,fill=indianred20197112,line width=0pt}
\addlegendentry{EWP}

\draw[draw=white,fill=indianred20197112,line width=0pt] (axis cs:0.75,0.00360291) rectangle (axis cs:1,0.01350303);
\draw[draw=white,fill=indianred20197112,line width=0pt] (axis cs:1.75,0.00538552) rectangle (axis cs:2,0.02022653);
\draw[draw=white,fill=indianred20197112,line width=0pt] (axis cs:2.75,0.0071511) rectangle (axis cs:3,0.02689982);
\draw[draw=white,fill=dimgray1419984,line width=0pt] (axis cs:-0.25,0.00677816) rectangle (axis cs:0,0.00986088);
\addlegendimage{ybar,ybar legend,draw=white,fill=dimgray1419984,line width=0pt}
\addlegendentry{L. Red}

\draw[draw=white,fill=dimgray1419984,line width=0pt] (axis cs:0.75,0.01350303) rectangle (axis cs:1,0.01814395);
\draw[draw=white,fill=dimgray1419984,line width=0pt] (axis cs:1.75,0.02022653) rectangle (axis cs:2,0.02645753);
\draw[draw=white,fill=dimgray1419984,line width=0pt] (axis cs:2.75,0.02689982) rectangle (axis cs:3,0.03467037);
\draw[draw=white,fill=rosybrown180124163,line width=0pt] (axis cs:-0.25,0.00986088) rectangle (axis cs:0,0.00988632);
\addlegendimage{ybar,ybar legend,draw=white,fill=rosybrown180124163,line width=0pt}
\addlegendentry{IFFT}

\draw[draw=white,fill=rosybrown180124163,line width=0pt] (axis cs:0.75,0.01814395) rectangle (axis cs:1,0.01817048);
\draw[draw=white,fill=rosybrown180124163,line width=0pt] (axis cs:1.75,0.02645753) rectangle (axis cs:2,0.02648361);
\draw[draw=white,fill=rosybrown180124163,line width=0pt] (axis cs:2.75,0.03467037) rectangle (axis cs:3,0.03469674);
\draw[draw=white,fill=darkkhaki213187103,line width=0pt] (axis cs:-0.25,0.00988632) rectangle (axis cs:0,0.00989687);
\addlegendimage{ybar,ybar legend,draw=white,fill=darkkhaki213187103,line width=0pt}
\addlegendentry{Unpad}

\draw[draw=white,fill=darkkhaki213187103,line width=0pt] (axis cs:0.75,0.01817048) rectangle (axis cs:1,0.01818114);
\draw[draw=white,fill=darkkhaki213187103,line width=0pt] (axis cs:1.75,0.02648361) rectangle (axis cs:2,0.02649418);
\draw[draw=white,fill=darkkhaki213187103,line width=0pt] (axis cs:2.75,0.03469674) rectangle (axis cs:3,0.03470759);
\draw[draw=white,fill=royalblue72120208,line width=0pt,postaction={pattern=north east lines, pattern color=white}] (axis cs:0,0) rectangle (axis cs:0.25,1.137e-05);
\draw[draw=white,fill=royalblue72120208,line width=0pt,postaction={pattern=north east lines, pattern color=white}] (axis cs:1,0) rectangle (axis cs:1.25,1.155e-05);
\draw[draw=white,fill=royalblue72120208,line width=0pt,postaction={pattern=north east lines, pattern color=white}] (axis cs:2,0) rectangle (axis cs:2.25,1.151e-05);
\draw[draw=white,fill=royalblue72120208,line width=0pt,postaction={pattern=north east lines, pattern color=white}] (axis cs:3,0) rectangle (axis cs:3.25,1.212e-05);
\draw[draw=white,fill=darkkhaki15817589,line width=0pt,postaction={pattern=north east lines, pattern color=white}] (axis cs:0,1.137e-05) rectangle (axis cs:0.25,3.305e-05);
\draw[draw=white,fill=darkkhaki15817589,line width=0pt,postaction={pattern=north east lines, pattern color=white}] (axis cs:1,1.155e-05) rectangle (axis cs:1.25,3.452e-05);
\draw[draw=white,fill=darkkhaki15817589,line width=0pt,postaction={pattern=north east lines, pattern color=white}] (axis cs:2,1.151e-05) rectangle (axis cs:2.25,3.38e-05);
\draw[draw=white,fill=darkkhaki15817589,line width=0pt,postaction={pattern=north east lines, pattern color=white}] (axis cs:3,1.212e-05) rectangle (axis cs:3.25,3.409e-05);
\draw[draw=white,fill=indianred20197112,line width=0pt,postaction={pattern=north east lines, pattern color=white}] (axis cs:0,3.305e-05) rectangle (axis cs:0.25,0.00331334);
\draw[draw=white,fill=indianred20197112,line width=0pt,postaction={pattern=north east lines, pattern color=white}] (axis cs:1,3.452e-05) rectangle (axis cs:1.25,0.00660463);
\draw[draw=white,fill=indianred20197112,line width=0pt,postaction={pattern=north east lines, pattern color=white}] (axis cs:2,3.38e-05) rectangle (axis cs:2.25,0.00988618);
\draw[draw=white,fill=indianred20197112,line width=0pt,postaction={pattern=north east lines, pattern color=white}] (axis cs:3,3.409e-05) rectangle (axis cs:3.25,0.01314469);
\draw[draw=white,fill=dimgray1419984,line width=0pt,postaction={pattern=north east lines, pattern color=white}] (axis cs:0,0.00331334) rectangle (axis cs:0.25,0.00637878);
\draw[draw=white,fill=dimgray1419984,line width=0pt,postaction={pattern=north east lines, pattern color=white}] (axis cs:1,0.00660463) rectangle (axis cs:1.25,0.01149697);
\draw[draw=white,fill=dimgray1419984,line width=0pt,postaction={pattern=north east lines, pattern color=white}] (axis cs:2,0.00988618) rectangle (axis cs:2.25,0.01658693);
\draw[draw=white,fill=dimgray1419984,line width=0pt,postaction={pattern=north east lines, pattern color=white}] (axis cs:3,0.01314469) rectangle (axis cs:3.25,0.0216393);
\draw[draw=white,fill=rosybrown180124163,line width=0pt,postaction={pattern=north east lines, pattern color=white}] (axis cs:0,0.00637878) rectangle (axis cs:0.25,0.0078116);
\draw[draw=white,fill=rosybrown180124163,line width=0pt,postaction={pattern=north east lines, pattern color=white}] (axis cs:1,0.01149697) rectangle (axis cs:1.25,0.01433727);
\draw[draw=white,fill=rosybrown180124163,line width=0pt,postaction={pattern=north east lines, pattern color=white}] (axis cs:2,0.01658693) rectangle (axis cs:2.25,0.02083983);
\draw[draw=white,fill=rosybrown180124163,line width=0pt,postaction={pattern=north east lines, pattern color=white}] (axis cs:3,0.0216393) rectangle (axis cs:3.25,0.02728756);
\draw[draw=white,fill=darkkhaki213187103,line width=0pt,postaction={pattern=north east lines, pattern color=white}] (axis cs:0,0.0078116) rectangle (axis cs:0.25,0.00805687);
\draw[draw=white,fill=darkkhaki213187103,line width=0pt,postaction={pattern=north east lines, pattern color=white}] (axis cs:1,0.01433727) rectangle (axis cs:1.25,0.01481609);
\draw[draw=white,fill=darkkhaki213187103,line width=0pt,postaction={pattern=north east lines, pattern color=white}] (axis cs:2,0.02083983) rectangle (axis cs:2.25,0.02155271);
\draw[draw=white,fill=darkkhaki213187103,line width=0pt,postaction={pattern=north east lines, pattern color=white}] (axis cs:3,0.02728756) rectangle (axis cs:3.25,0.02823265);
\draw[draw=white,fill=gray,very thin] (axis cs:-0.4,0) rectangle (axis cs:0.4,0);
\draw[draw=white,fill=gray,very thin,postaction={pattern=north east lines, pattern color=white}] (axis cs:-0.4,0) rectangle (axis cs:0.4,0);
  \end{axis}
  \node [fill opacity=0.8,
    draw opacity=1,
    text opacity=1,
    draw=lightgray204,
    fill=gainsboro229,
    yshift=-1cm,
    anchor=west,name=leg2] at (leg.south west) {\shortstack[l]{
    \adjustbox{cframe=gray, bgcolor=gray}{\textcolor{gray}{\bf{\tiny{///}}}} $\blocktoep\,\,\,$ Matvec \\
\adjustbox{cframe=gray, bgcolor=gray}{\textcolor{white}{\bf{\tiny{///}}}} $\blocktoep^*$ Matvec }};
\end{tikzpicture}